\documentclass[11pt]{article}

\usepackage{amsfonts}
\usepackage{xcolor}
\usepackage{amsmath,amssymb,comment}
\usepackage{verbatim}

\usepackage{tikz}
\usetikzlibrary{shapes,arrows}
\usetikzlibrary{intersections}
\usepackage{pgfplots}
\pgfplotsset{compat = newest}
\newtheorem{theo}{Theorem}[section]
\newtheorem{lem}[theo]{Lemma}

\newtheorem{prop}[theo]{Proposition}

\newcommand{\mysection}[1]{\section{#1} \setcounter{equation}{0}}
    \makeatletter
\def\@fnsymbol#1{\ensuremath{\ifcase#1\or *\or \ddagger\or
   \mathsection\or \mathparagraph\or \|\or **\or \dagger\dagger
   \or \ddagger\ddagger \else\@ctrerr\fi}}
    \makeatother
\newcommand{\proof}{{\sc Proof.} \quad}
\newcommand{\proofc}{{\sc Proof} \ }
\newcommand{\be}{\begin{equation} \label}
\newcommand{\ee}{\end{equation}}
\newcommand{\bea}{\begin{eqnarray}\label}
\newcommand{\eea}{\end{eqnarray}}
\newcommand{\bas}{\begin{eqnarray*}}
\newcommand{\eas}{\end{eqnarray*}}
\newcommand{\bit}{\begin{itemize}}
\newcommand{\eit}{\end{itemize}}
\newcommand{\qed}{\hfill$\Box$ \vskip.2cm}
\newcommand{\nn}{\nonumber}
\newcommand{\R}{\mathbb{R}}

\newcommand{\N}{\mathbb{N}}

\newcommand{\eps}{\varepsilon}

\newcommand{\supp}{{\rm supp} \, }

\newcommand{\na}{\nabla}
\newcommand{\Del}{\Delta}
\newcommand{\del}{\delta}
\newcommand{\al}{\alpha}
\newcommand{\lam}{\lambda}

\newcommand{\om}{\omega}

\newcommand{\Mint}{- \hspace*{-4mm} \int}

\newcommand{\ov}{\overline}
\newcommand{\un}{\underline}
\newcommand{\wh}{\widehat}

\newcommand{\hs}{\hspace*}
\newcommand{\sm}{\setminus}

\newcommand{\vp}{\varphi}
\newcommand{\vt}{\vartheta}
\newcommand{\lbal}{\left\{ \begin{array}{l}}
\newcommand{\lball}{\left\{ \begin{array}{ll}}
\newcommand{\ear}{\end{array} \right.}

\newcommand{\cred}{\normalcolor}

\newcommand{\abs}{\\[5pt]}
\newcommand{\Abs}{\\[5mm]}

\newcommand{\ow}{\ov{w}}

\newcommand{\zeps}{z_\eps}

\newcommand{\ze}{\zeta_\eps}
\newcommand{\zes}{\zeta_{\eps s}}
\newcommand{\zess}{\zeta_{\eps ss}}
\newcommand{\K}{K_\chi}
\newcommand{\ophi}{\ov{\vp}}
\newcommand{\uphi}{\un{\vp}}
\newcommand{\xe}{\xi_\eps}
\newcommand{\xes}{\xi_{\eps s}}
\newcommand{\xess}{\xi_{\eps ss}}
\newcommand{\xet}{\xi_{\eps t}}
\newcommand{\whs}{\wh{s}}
\newcommand{\wht}{\wh{t}}
\begin{document}
\enlargethispage{10mm}
\title{Stability vs.~instability of singular steady states\\
in the parabolic-elliptic Keller-Segel system on $\R^n$}
\author{
Francesca Colasuonno\footnote{francesca.colasuonno@unibo.it}\\
{\small Dipartimento di Matematica, Alma Mater Studiorum Universit\`a di Bologna,}\\
{\small Piazza di Porta San Donato 5, 40126 Bologna, Italy}
\and
Michael Winkler\footnote{michael.winkler@math.uni-paderborn.de}\\
{\small Institut f\"ur Mathematik, Universit\"at Paderborn,}\\
{\small 33098 Paderborn, Germany} 
}
\date{}
\maketitle
\begin{abstract}
\noindent 
The Cauchy problem in $\R^n$ is considered for 
\bas
	\lbal
	u_t = \Del u - \na \cdot (u\na v), \\[1mm]
	0 = \Del v + u.
	\ear
\eas
For each $n\ge 10$, a statement on stability and attractiveness of the singular steady state given by
\bas
	u_\star(x):=\frac{2(n-2)}{|x|^2},
	\qquad x\in\R^n\sm\{0\},
\eas
is derived within classes of nonnegative radial solutions emanating from initial data less concentrated than $u_\star$.
In particular, for any such $n$ it is shown that infinite-time blow-up occurs for 
all radial initial data which are less concentrated than $u_\star$ and satisfy
\bas
	u_0(x) \ge \frac{2(n-2)}{|x|^2} - \frac{C}{|x|^{2+\theta}}
	\qquad \mbox{for all } x\in \R^n\sm B_1(0)
\eas
with some $C>0$ and some $\theta>\frac{n-2+\sqrt{(n-2)(n-10)}}{2}$.\abs
This is complemented by a result which, in the case when $3\le n \le 9$, asserts 
instability of $u_\star$ as well as the existence of a bounded absorbing set for
all radial trajectories initially less concentrated than $u_\star$.\abs
{\cred
In particular, previous knowledge on stability properties of $u_\star$, as having been gained for $n\ge 11$ in \cite{senba_FE2013},
is thereby extended to any dimension $n\ge 3$.
}
\abs
\end{abstract}
\newpage
\section{Introduction}\label{intro}
This manuscript studies the parabolic-elliptic Keller-Segel system 
\be{0}
	\left\{ \begin{array}{ll}	
	u_t = \Del u - \na\cdot (u\na v), 
	\qquad & x\in\R^n, \ t>0, \\[1mm]
	0 = \Del v + u,
	\qquad & x\in\R^n, \ t>0, \\[1mm]
	u(x,0)=u_0(x),
	\qquad & x\in\R^n,
	\end{array} \right.
\ee
where $n\ge 3$ and $u_0\in C^0(\R^n)$ is nonnegative and bounded.
Systems of this form
play a key role not only in simplified mathematical descriptions of bacterial aggregation, but also in the modeling 
of gravitational collapse in stellar dynamics \cite{KS,JL,chandra42,biler_book,biler1995}.\abs
Accordingly, a predominant focus in analytical studies concerned with (\ref{0}) is 
on aspects related to the singularity-supporting potential of the interaction mechanism therein. 
Indeed, the literature has provided a meanwhile rich collection of results which identify constellations under which
solutions to (\ref{0}) either remain globally bounded,
or reflect finite-time formation of singularities in the sense that 
$\limsup_{t\nearrow T} \|u(\cdot,t)\|_{L^\infty(\R^n)}=\infty$ for some $T\in (0,\infty)$,
see \cite{biler_karch_zienkiewicz_NHM,biler_karch_pilarczyk,senba_FE2005,perthame_EJDE,naito_senba,naito_JDE};
cf.~also the broader discussions in \cite{biler_boritchev_brandolese,suzuki_book,biler_book}.\abs
A natural next step consists in the investigation of qualitative solution behavior near singularities.
In this regard, considerably precise characterizations are available in cases when a special ansatz, e.g.~of self-similar type,
is pursued \cite{herrero_medina_velazquez,naito_suzuki,senba_FE2005,collot_et_al};
apart from this, quite comprehensive knowledge on temporal rates of blow-up has been achieved 
\cite{mizoguchi_senba2007,senba_ADE2006,naito_senba}. In particular, it is possible to prove that solutions $(u,v)$ which blow up at finite time $T$ satisfy the following lower bound: $\|u(\cdot,t)\|_{L^\infty(\mathbb R^n)}\ge (T-t)^{-1}$. Then, the blow-up is called of type I if $\limsup_{t\nearrow T} (T-t)\|u(\cdot,t)\|_{L^\infty(\R^n)}< \infty$, otherwise it is called of type II. The dynamics of the system and its analysis strongly depend on the dimension $n$. In dimension one, solutions exist globally in time, so no finite time blow-up is allowed, while for dimensions $n\ge 2$ different behaviors are possible. In dimension two, the situation changes drastically depending on whether the mass of the system is smaller, equal, or larger than the critical threshold $8\pi$: in the first case solutions are global in time, in the latter they blow up in finite time, while, when the mass is critical, infinite time blow-up may occur. 
For critical mass phenomena in bounded domains $\Omega\subset\mathbb R^n$, we refer to \cite{JL} for the planar case $n=2$, and to the recent paper \cite{win_MATANN} for the radial case $\Omega=B_R(0)$, in possibly higher dimensions $n\ge 2$. 
Radial solutions of \eqref{0} exhibiting type II blow-up are found in \cite{HerreroVelazquez} for $n=2$ and in \cite{senba_ADE2006} for $n\ge 11$. When $3\le n\le 9$, sufficient conditions for type I blow-up are given in \cite{mizoguchi_senba2007}; cf.~also \cite{brenner}.\abs
A distinctly simple example of singular behavior in (\ref{0})
is formed by the explicit steady state solutions, also referred to as the Chandrasekhar solutions \cite{biler_book}, 
which are given by		
\be{stat}
	u(x,t)=u_\star(x)
	:= \frac{2(n-2)}{|x|^2},
	\qquad x\in\R^n\sm\{0\}, \ t>0.
\ee
They do not only hint, through their mere presence as nontrivial equlibria when $n\ge 3$, at several fundamental differences 
between the two- and the higher-dimensional versions of (\ref{0}) (see \cite{perthame_EJDE,senba_FE2013}, for instance),
but beyond this the particular functional form of the balance between diffusion and cross-diffusion expressed in \eqref{stat}
seems to quantify key properties of the separatrix distinguishing between attraction-driven and diffusion-dominated behavior.
For instance, the decay asymptotics in (\ref{stat}) appear to mark 
a critical tail behavior with respect to the occurrence of explosions 
\cite{naito_JDE,biler_karch_zienkiewicz_NHM,win_slowly_decaying_data}; cf.~also Proposition \ref{prop0} below.

It is worth mentioning that the Chandrasekhar solution $u_\star$ plays an important role also for the problem set in bounded domains. Indeed, recalling that, given two nonnegative, radial functions $\varphi$ and $\psi\in L^1(B_R(0))$, $\varphi$ is said to be less concentrated than $\psi$ if $\int_{B_r(0)}\varphi\le \int_{B_r(0)}\psi$ for every $r\in(0,R)$, in \cite[Theorem 1.3]{win_MATANN} it is shown that the problem set in $B_R(0)$ admits a global bounded classical solution provided that the initial datum $u_0\in C^0(\overline{B_R}(0))$ is nonnegative, radial, and less concentrated than $u_\star$.\abs
{\bf The goal: Characterizing stability properties of $u_\star$.} \quad
The present manuscript now addresses the question how far the structure of the stationary point $u_\star$
is of perceptible relevance for the behavior of solutions emanating from nearby initial data. 
Mainly due to the singular nature of $u_\star$, this issue appears to go beyond the scope of straightforward
stability analysis based on standard linearization, especially since any constellation involving nontrivial
domains of attraction should bring about some solutions to (\ref{0}) which undergo a genuine explosion.
After all, when focusing on deviations from $u_\star$ which in a suitable sense lead to less concentrated distributions,
we may rely on the following essentially well-known basic statement on global classical solvability
to exclude the possibility of finite-time blow-up
(see, e.g., \cite{naito_JDE} and also \cite[Proposition 1.1, Theorem 1.2]{win_slowly_decaying_data}).
\begin{prop}\label{prop0}
  Let $n\ge 3$, and assume that
  \be{init}
	u_0\in C^0(\R^n) \mbox{\; is radially symmetric and such that \; }
	0 \le u_0(x) \le \frac{2(n-2)}{|x|^2}\;
	\mbox{ for all } x\in\R^n\sm \{0\}.
  \ee
  Then there exist radial functions
  \bas
	\lbal
	u\in C^0(\R^n\times [0,\infty)) \cap C^\infty(\R^n\times (0,\infty))
	\qquad \mbox{and} \\[1mm]
	v\in C^\infty(\R^n\times (0,\infty)),
	\ear
  \eas
  with $u$ being uniquely determined by the additional requirements that $u\ge 0$ and that
  \bas
	\sup_{t\in (0,T)} \|u(\cdot,t)\|_{L^\infty(\R^n)} < \infty
	\qquad \mbox{for all } T>0,
  \eas
  such that (\ref{0}) is satisfied in the classical sense.
  Moreover,
  \be{0.1}
	\Mint_{B_r(0)} u(x,t) dx \le \frac{2n}{r^2}
	\qquad \mbox{for all $r>0$ and } t>0.
  \ee
\end{prop}
The latter leaves open the possibility of explosions to occur in infinite time, and in fact it has been derived
in \cite[Corollary 1.2]{senba_FE2013} that when $n\ge 11$, some smooth initial data fulfilling \eqref{init} lead to solutions
for which, although $\liminf_{t\to\infty} \|u(\cdot,t)\|_{L^\infty(\R^n)} < \infty$, we have
$\limsup_{t\to\infty} \|u(\cdot,t)\|_{L^\infty(\R^n)} =\infty$
(cf.~also \cite{senba_NA2009} for a different construction of global unbounded solutions).\Abs
{\bf Main results.} \quad
The first of our main results now reveals that actually for each $n\ge 10$, the singular states in (\ref{stat}) indeed attract
considerably many solutions. Through its nature as a statement on asymptotic stability, 
this does not only assert the occurrence of global solutions which grow up in an essentially
non-oscillatory manner described by \eqref{7.3}, but through the assertion on convergence in \eqref{7.4}
this moreover also provides information on the asymptotic profile. 
\begin{theo}\label{theo7}
  Let $n\ge 10$, let $\theta>0$ be such that
  \be{theta}
	\theta>\frac{n-2+\sqrt{(n-2)(n-10)}}{2},
  \ee
  and suppose that $u_0$ satisfies (\ref{init}) as well as
  \be{7.2}
	u_0(x) \ge \frac{2(n-2)}{|x|^2} - \frac{C}{|x|^{2+\theta}}
	\qquad \mbox{for all } x\in \R^n\sm B_1(0)
  \ee
  with some $C>0$.
  Then the solution $(u,v)$ of (\ref{0}) from Proposition \ref{prop0} blows up in infinite time in the sense that
  \be{7.3}
	\|u(\cdot,t)\|_{L^\infty(\R^n)} \to \infty
	\qquad \mbox{as } t\to\infty.
  \ee
  Moreover, with $u_\star$ taken from (\ref{stat}) we have
  \be{7.4}
	u(\cdot,t) \to u_\star
	\quad \mbox{in } C^0_{loc}(\R^n\setminus\{0\})
	\qquad \mbox{as } t\to\infty.
  \ee
\end{theo}
We remark here that under the slightly stronger assumption $n\ge 11$, 
for each member of a family $(u_\kappa)_{\kappa\ge 0}$ of bounded steady state solutions to (\ref{0}) satisfying
$u_\kappa(0)=\kappa$, 
a corresponding stability analysis in the respective setting of regular functions can be built on parabolic comparison 
with suitable bounded sub- and supersolutions; 
an accordingly resulting statement on attractiveness comparable to that from Theorem \ref{theo7}
can be found documented in \cite{senba_FE2013}.\abs
Now in stark contrast to the situation in Theorem \ref{theo7}, in any space of lower dimension the equilibrium $u_\star$
turns out to become quite strongly unstable with respect to the class of radial and concentration-reducing perturbations
considered above.
Most drastically, this repulsive character becomes manifest in the following statement which does not only rule 
out any unboundedness phenomenon in such cases, but moreover also asserts the existence of a bounded absorbing set 
for all solutions to (\ref{0}) within the range of initial data which are less concentrated than $u_\star$ by 
satisfying (\ref{init}).
\begin{theo}\label{theo25}
  Let $3\le n\le 9$. 
  Then one can find $C>0$ with the property that whenever $u_0$ satisfies (\ref{init}), there exists $t_\star=t_\star(u_0)>0$
  such that the global classical solution $(u,v)$ of (\ref{0}) from Proposition \ref{prop0} satisfies
  \[
	\|u(\cdot,t)\|_{L^\infty(\R^n)} \le C
	\qquad \mbox{for all } t>t_\star.
  \]
%
\end{theo}
Apart from that, however, also far outside the spatial origin such solutions will eventually be repelled to a significant
extent:
\begin{prop}\label{prop26}
  If $3\le n \le 9$, then for all $r_0>0$ there exists $C(r_0)>0$ such that if (\ref{init}) holds, then 
  it is possible to find $t_{\star\star}=t_{\star\star}(u_0)>0$ satisfying
  \[
	\Mint_{B_r(0)} u_\star(x) dx - \Mint_{B_r(0)} u(x,t) dx
	\ge C(r_0)
	\qquad \mbox{for all $r\in (0,r_0)$ and } t> t_{\star\star}.
  \]
\end{prop}
{\cred
In summary, within the realm of radially symmetric solutions 
the above results seem to provide a fairly complete picture with respect to stability and attractiveness of the
singular steady state $u_\star$ in all dimensions $n\ge 3$.
An interesting issue left unaddressed here is how far comparable phenomena can be observed also when nonradial initial
data are included. 
}

\bigskip
{\bf Main ideas.} \quad
In \cite{senba_FE2013}, the mainly comparison-based stability analysis of regular equilibria crucially relied
on the assumption $n\ge 11$ through a requirement
on positivity of the radicand in \eqref{theta}, and hence on distinctness of two roots of a related quadratic polynomial.
In order to design an approach capable of covering also the resonance case $n=10$ in which both these roots collapse,
the strategy toward our verification of Theorem \ref{theo7} will be based on a variational idea.
At its core, namely, our analysis in this first manuscript part in Section \ref{sect3}
aims at deriving a rigorous counterpart of the energy inequality
\be{energy}
	\frac{1}{p} \phi'(t) 
	\le - \frac{n\gamma}{p+1} \int_0^\infty s^{-\gamma-1} \vp^{p+1}(s,t) ds,
	\qquad t>0,
\ee
which for $n\ge 10$, and within suitable ranges of the auxiliary parameters $p$ and $\gamma$, 
is formally associated with (\ref{0}) along trajetories fulfilling (\ref{init}) (see Lemma \ref{lem44}), where
\be{funct}
	\phi(t):=\int_0^\infty s^{-\gamma} \vp^p(s,t) ds,
	\qquad t\ge 0,
\ee
with 
\bas
	\vp(s,t):=\int_0^{s^\frac{1}{n}} \rho^{n-1} \big\{ u_\star(\rho) - u(\rho,t)\big\} d\rho,
	\qquad s\ge 0, t\ge 0,
\eas
representing the deviation of the radial solution $u=u(r,t)$ from $u_\star$, formulated here in terms of cumulated densities.
Here, substantial use will be made of the exhaustive knowledge on the precise functional form of the equilibrium $u_\star$,
according to which the function $\vp$, nonnegative due to a comparison argument, is known to solve a scalar parabolic
problem with fully explicit ingredients (see \eqref{0p}). 
In a first stage of our reasoning, we may then merely rely on the dissipation mechanism expressed in \eqref{energy} to conclude
that for appropriately chosen $p$ and $\gamma$ we have
\bas
	\int_0^\infty \int_0^\infty s^{-\gamma-1} \vp^{p+1}(s,t) dsdt < \infty
\eas
(Lemma \ref{lem44}),
whereupon only in a second step, by means of additional regularity properties of $\vp$ due to parabolic smoothing,
the full topological setting from Theorem \ref{theo7} will be addressed (Lemma \ref{lem45}). \abs
The approach underlying our derivation of the instability results in Theorem \ref{theo25} and Proposition \ref{prop26} will
be entirely different by being exclusively based on comparison methods in all essential parts:
In Section \ref{sect4}, a first parabolic comparison, crucially relying on the hypothesis $n\le 9$,
involves separated subsolutions containing compactly supported  cutouts of spatially oscillating functions in their spatial part, and certain solutions of logistic-type Bernoulli ODEs; 
a basic instability property will here result from the repulsive character of its trivial equilibrium (Lemma \ref{lem21}).
The spatially local information thereby gained is extended to balls of arbitrary fixed size 
by means of a second comparison argument in Lemma \ref{lem22},
and an upshift of corresponding estimates, transferring bounds 
from cumulated quantities to their derivatives and hence to the original densities,
is finally achieved by means of a third application of a maximum principle in the course of a Bernstein-type argument (claim \eqref{24.1} in the proof of Theorem~\ref{theo25}).
%
%
%
%
%
%
%
%
%
%
%
%
%
%
%
%
\mysection{Some preparations}
Let us first record some elementary features of the mass accumulation functions associated with solutions of (\ref{0}),
particularly with regard to their deviation from the singular state in (\ref{stat}).
In standard abuse of notation, here and throughout 	
the sequel we shall switch to radial notation when writing, e.g., $u(r,t)$ instead of $u(x,t)$ for $r=|x|\ge 0$.
\begin{lem}\label{lem2}
  Let $n\ge 3$, assume (\ref{init}), and let $(u,v)$ be as in Proposition \ref{prop0}.
  Then letting
  \be{w}
	w(s,t):=\int_0^{s^\frac{1}{n}} \rho^{n-1} u(\rho,t) d\rho, 
	\qquad s\ge 0, t\ge 0,
  \ee
  and
  \be{phi}
	\vp(s,t):=w_\star(s)-w(s,t),
	\qquad s\ge 0, t\ge 0,
  \ee
  with
	\[	
	w_\star(s):=2s^{1-\frac{2}{n}},
	\qquad s\ge 0,
	\] 
  defines functions $w$ and $\vp$ which belong to $C^0([0,\infty)^2) \cap C^{2,1}((0,\infty)^2)$ and satisfy
  \be{0w}
	w_t = n^2 s^{2-\frac{2}{n}} w_{ss} + nww_s,
	\qquad s>0, \ t>0,
  \ee
  and
  \be{0p}
	\vp_t = n^2 s^{2-\frac{2}{n}} \vp_{ss} + 2n s^{1-\frac{2}{n}} \vp_s + 2(n-2) s^{-\frac{2}{n}} \vp - n\vp \vp_s,
	\qquad s>0, \ t>0,
  \ee
  as well as
  \be{pb}
	0 \le w(s,t) \le w_\star(s)
	\quad \mbox{and} \quad
	0 \le \vp(s,t) \le w_\star(s)
	\qquad \mbox{for all $s\ge 0$ and } t\ge 0.
  \ee
\end{lem}
\proof
  This can be seen by straightforward computation using the radial version of (\ref{0}). Indeed, performing the partial derivative in $s$ of the parameter-dependent integral that defines $w$, we 
{\cred see that
}
  \[w_s(s,t)=\rho^{n-1}u(\rho,t)\Big|_{\rho=s^{\frac{1}{n}}}(s^{\frac{1}{n}})'=\frac{1}{n}u(s^{\frac{1}{n}},t),\]
  whence
  \[w_{ss}(s,t)=\frac{1}{n^2}s^{\frac{1}{n}-1}u_r(s^{\frac{1}{n}},t).\]
  Reasoning as in \cite[Section 4]{JL}, we integrate the first equation of \eqref{0} over the ball $B_{s^{\frac{1}{n}}}(0)$ to get
  \[\omega_{n-1}\int_0^{s^{\frac{1}{n}}}u_t(r,t)r^{n-1}dr=\omega_{n-1}s^{\frac{n-1}{n}}u_r(s^{\frac{1}{n}},t)-\omega_{n-1}u(s^{\frac{1}{n}},t)v_r(s^{\frac{1}{n}},t),\] 
  where $\omega_{n-1}$ denotes the $(n-1)$-dimensional measure of the unit ball in $\mathbb R^n$ and we applied twice the the Divergence Theorem.
Thus, keeping in mind the definition of $w$ and the expressions of $w_s$ and of $w_{ss}$ written above, and using the second equation of \eqref{0} integrated over $B_{s^{\frac{1}{n}}}(0)$, we obtain \eqref{0w}.
As for \eqref{0p}, by \eqref{0w} we can compute 
\[
\begin{aligned}
\vp_t& =  -n^2s^{2-\frac{2}{n}}w_{ss}-nww_s = n^2 s^{2-\frac{2}{n}}\vp_{ss}-n\vp \vp_s
+2ns^{1-\frac{2}{n}}\vp_s+2(n-2)s^{-\frac{2}{n}}\vp,\end{aligned}
\]
where we expressed $w_{ss}$ in terms of $\vp_{ss}$ using $\vp_{ss}=(w_\star)_{ss}-w_{ss}=-\frac{4}{n}\left(1-\frac{2}{n}\right)s^{-1-\frac{2}{n}}-w_{ss}$. This concludes the proof of \eqref{0p}. Finally, translating \eqref{0.1}
{\cred yields 
}
\eqref{pb}.
\qed
To prepare a unified treatment of several cut-off procedures below, let us fix $\chi\in C^\infty([0,\infty))$ such that
\be{chi}
	0 \le \chi \le 1
	\mbox{ and } \chi'\ge 0
	\mbox{ on } [0,\infty),
	\quad
	\chi\equiv 0
	\mbox{ on } [0,1]
	\quad \mbox{and} \quad
	\chi\equiv 1
	\mbox{ on } [2,\infty),
\ee
and set
\be{ze}
	\ze(s):=\lball
	\chi(\frac{2s}{\eps}),
	\qquad s\in [0,1], \\[1mm]
	1-\chi(\eps s),
	\qquad s>1,
	\ear
\ee
for $\eps\in (0,1)$.
\begin{center}
 \pgfplotsset{width=7cm}
\begin{tikzpicture}[
declare function={
    func(\x)= (\x < 1) * 0   +
              and(\x >= 1, \x<2) * e^(1+1/((\x-2)^2-1))     +
              (\x>= 2) * 1 
   ;
    gunc(\x)= and(\x >=0, \x< 1) * func(2*\x/0.9)   +
              (\x>= 1) * (1 - func(\x*0.9)
   ;
  } 
  ]
\begin{axis}[
  axis equal,
  grid=major,
  axis x line=center,
  axis y line=center,
  xtick={1},
   extra x ticks={0.45,0.9,1.18,2.36},
    extra x tick style={grid=major},
    extra x tick labels={$\frac{\varepsilon}{2}$,$\varepsilon$,$\frac{1}{\varepsilon}$,$\frac{2}{\varepsilon}$},
  ytick={1},
  xlabel={$s$}, xlabel style={below right},
  ylabel={$\zeta_\varepsilon$}, 
  xmin=0,
  xmax=2.5,
  ymin=0,
  ymax=1.0,
%
 ]
\addplot[black,thick,domain=0:2.5,samples=1000]{gunc(x)};
\end{axis}
\end{tikzpicture}

{\it
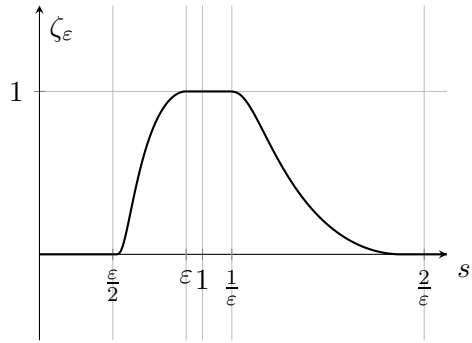
Fig.~1: Shape of the cut-off function $\zeta_\eps$
}
\end{center}

Then $(\ze)_{\eps\in (0,1)} \subset C^\infty([0,\infty))$ with
\be{z2}
	0 \le \ze \le 1
	\mbox{ on } [0,\infty),
	\quad
	\ze\equiv 1
	\mbox{ on $[\eps,\frac{1}{\eps}]$}
	\quad \mbox{and} \quad
	\supp \ze \subset \mbox{$[\frac{\eps}{2},\frac{2}{\eps}]$}
	\qquad \mbox{for all } \eps\in (0,1),
\ee
and writing 
\be{K}
	\K:=\max \Big\{ 2\|\chi'\|_{L^\infty((0,\infty))} \, , \, 4\|\chi''\|_{L^\infty((0,\infty))}\Big\}, 
\ee
we obtain that besides
\be{chi1}
	|\chi'(t)|\le \K
	\qquad \mbox{for all } t>0,
\ee
we also have
\be{z3}
	|\zes(s)| \le \frac{\K}{\eps}
	\quad \mbox{and} \quad
	|\zess(s)|\le\frac{\K}{\eps^2}
	\qquad \mbox{for all $s \in (\frac{\eps}{2},\eps)$ and } \eps\in (0,1)
\ee
as well as
\be{z4}
	|\zes(s)| \le \K \eps
	\quad \mbox{and} \quad
	|\zess(s)|\le \K \eps^2
	\qquad \mbox{for all $s\in (\frac{1}{\eps},\frac{2}{\eps})$ and } \eps\in (0,1).
\ee
{\cred
The role and this particular choice of the cut-off functions $\ze$ will become apparent especially in 
}
the proofs of Lemmas \ref{lem42} and \ref{lem43}, and of Theorem \ref{theo25}.

\mysection{Stability and attractiveness of $u_\star$. Proof of Theorem \ref{theo7}}\label{sect3}
{\cred
The intention of this section is to derive the stability result in Theorem \ref{theo7} by means of suitable dissipative
properties of the singularly weighted functional in \eqref{funct}.
Our considerations in this direction will be launched by the following observation that
describes a basic evolution feature enjoyed by a regularized variant thereof.
}
\begin{lem}\label{lem4}
  Let $n\ge 3$ and assume (\ref{init}),
  and let $p>1$ and $\gamma>0$. Then for each $\eps\in (0,1)$, letting
  \be{pe}
	\phi_\eps(t):=\int_0^\infty s^{-\gamma} \ze^2(s) \vp^p(s,t) ds,
	\qquad t\ge 0,
  \ee
  with $\vp$ and $\ze$ as in \eqref{phi} and \eqref{ze}, we obtain a nonnegative function $\phi_\eps \in C^0([0,\infty)) \cap
  C^1((0,\infty))$ which satisfies
  \bea{4.1}
	\frac{1}{p}\phi_\eps'(t)
	&=& - n^2(p-1) \int_0^\infty s^{2-\frac{2}{n}-\gamma} \ze^2 \vp^{p-2} \vp_s^2 ds \nn\\
	& & + \Big\{ \frac{n}{p} (n\gamma-2n+4)\Big(\gamma-1+\frac{2}{n}\Big) + 2(n-2) \Big\} \cdot 
		\int_0^\infty s^{-\frac{2}{n}-\gamma} \ze^2 \vp^p ds \nn\\
	& & - \frac{n\gamma}{p+1} \int_0^\infty s^{-\gamma-1} \ze^2 \vp^{p+1} ds \nn\\
	& & - \frac{4n}{p} (n\gamma-2n+3) \int_0^\infty s^{1-\frac{2}{n}-\gamma} \ze \zes \vp^p ds \nn\\
	& & + \frac{2n^2}{p} \int_0^\infty s^{2-\frac{2}{n}-\gamma} (\ze \zess + \zes^2) \vp^p ds \nn\\
	& & + \frac{2n}{p+1} \int_0^\infty s^{-\gamma} \ze \zes \vp^{p+1} ds
	\qquad \mbox{for all } t>0.
  \eea
\end{lem}
\proof
{\cred
  Since $(0,\infty) \ni s \mapsto s^{-\gamma} \ze^2(s)$ 
  is a time-independent smooth function with
}
  compact support in $(0,\infty)$, 
  the claimed regulatity features immediately result from 
  Lemma \ref{lem2}. 
  To derive \eqref{4.1}, we only need to use \eqref{0p} in computing
  \bea{4.2}
	\frac{1}{p} \phi_\eps'(t)
	&=& \int_0^\infty s^{-\gamma} \ze^2 \vp^{p-1} \vp_t ds \nn\\
	&=& n^2 \int_0^\infty s^{2-\frac{2}{n}-\gamma} \ze^2 \vp^{p-1} \vp_{ss} ds
	+ 2n \int_0^\infty s^{1-\frac{2}{n}-\gamma} \ze^2 \vp^{p-1} \vp_s ds \nn\\
	& & + 2(n-2) \int_0^\infty s^{-\frac{2}{n}-\gamma} \ze^2 \vp^p ds
	- n \int_0^\infty s^{-\gamma} \ze^2 \vp^p \vp_s ds 
	\qquad \mbox{for all } t>0,
  \eea
  and to integrate by parts in verifying that
  \bas
	n^2 \int_0^\infty s^{2-\frac{2}{n}-\gamma} \ze^2 \vp^{p-1} \vp_{ss} ds
	&=& - n^2(p-1) \int_0^\infty s^{2-\frac{2}{n}-\gamma} \ze^2 \vp^{p-2} \vp_s^2 ds \nn\\
	& & - n^2\Big(2-\frac{2}{n}-\gamma\Big) \int_0^\infty s^{1-\frac{2}{n}-\gamma} \ze^2 \vp^{p-1} \vp_s ds \nn\\
	& & - 2n^2 \int_0^\infty s^{2-\frac{2}{n}-\gamma} \ze \zes \vp^{p-1} \vp_s ds
	\qquad \mbox{for all } t>0,
  \eas
  that
  \bea{4.3}
	& & \hs{-30mm}
	\Big\{ - n^2\Big(2-\frac{2}{n}-\gamma\Big) + 2n\Big\} \cdot
	\int_0^\infty s^{1-\frac{2}{n}-\gamma} \vp^{p-1} \vp_s ds \nn\\
	&=& -\frac{n}{p} (2n-4-n\gamma) \int_0^\infty  s^{1-\frac{2}{n}-\gamma} \ze^2 (\vp^p)_s ds \nn\\
	&=& \frac{n}{p} (2n-4-n\gamma) \Big(1-\frac{2}{n}-\gamma\Big) \int_0^\infty s^{-\frac{2}{n}-\gamma} \ze^2 \vp^p ds \nn\\
	& & + \frac{2n}{p} (2n-4-n\gamma) \int_0^\infty s^{1-\frac{2}{n}-\gamma} \ze \zes \vp^p ds
	\qquad \mbox{for all } t>0
  \eea
  and, similarly,
  \bea{4.4}
	\hs{-6mm}
	-2n^2 \int_0^\infty s^{2-\frac{2}{n}-\gamma} \ze \zes \vp^{p-1} \vp_s ds
	&=& \frac{2n^2}{p}\Big(2-\frac{2}{n}-\gamma\Big) \int_0^\infty s^{1-\frac{2}{n}-\gamma} \ze\zes \vp^p ds \nn\\
	& & + \frac{2n^2}{p} \int_0^\infty s^{2-\frac{2}{n}-\gamma} (\ze\zess+\zes^2) \vp^p ds 
	\qquad \mbox{for all } t>0
  \eea
  as well as
  \bas	
	- n \int_0^\infty s^{-\gamma} \ze^2 \vp^p \vp_s ds
	&=& - \frac{n\gamma}{p+1} \int_0^\infty s^{-\gamma-1} \ze^2 \vp^{p+1} ds \nn\\
	& & + \frac{2n}{p+1} \int_0^\infty s^{-\gamma} \ze\zes \vp^{p+1} ds
	\qquad \mbox{for all } t>0.
  \eas
  Linking the second last summands in \eqref{4.2} and \eqref{4.3}, and combining the last expression in (\ref{4.3}) with the second 
  last from \eqref{4.4}, namely, readily shows that since
  \bas
	\frac{2n}{p}(2n-4-n\gamma) + \frac{2n^2}{p}\Big(2-\frac{2}{n}-\gamma\Big)
	= \frac{2n}{p}(2n-4-n\gamma+2n-2-n\gamma)
	= - \frac{4n}{p}(n\gamma-2n+3),
  \eas
  the identity in \eqref{4.1} is equivalent to (\ref{4.2}).
\qed
Now, up to expressions containing derivatives of the cut-off functions $\ze$, in the framework of Theorem~\ref{theo7} and
within a suitable range of $p$, the second summand on the right of (\ref{4.1})
can be estimated against the diffusion-related first one on the basis of the Hardy inequality.
This crucial step, to be completed in Lemma \ref{lem44}, will be prepared by the following.
\begin{lem}\label{lem41}
  Suppose that $n\ge 3$ and that (\ref{init}) holds, and let $p>1$ and $\gamma>1-\frac{2}{n}$.
  Then with $(\ze)_{\eps\in (0,1)}$ and $\vp$ as in \eqref{ze} and \eqref{phi},
  \bea{41.1}
	\int_0^\infty s^{-\frac{2}{n}-\gamma} \ze^2 \vp^p ds
	&\le\displaystyle{\frac{1}{(\gamma-1+\frac{2}{n})^2}}&\left\{p^2\int_0^\infty s^{2-\frac{2}{n}-\gamma} \ze^2 \vp^{p-2} \vp_s^2 ds \right. \nn\\
	& &\quad +4(\gamma-2+\frac{2}{n}) \int_0^\infty s^{1-\frac{2}{n}-\gamma} \ze\zes \vp^p ds \nn\\
	& &\left.\quad - 4 \int_0^\infty s^{2-\frac{2}{n}-\gamma} \ze \zess \vp^p ds\right\}
  \eea
  for all $t>0$ and $\eps\in (0,1)$.
\end{lem}
\proof
  An application of the Hardy inequality,
  \bas
	\int_0^\infty s^{-\beta} \psi^2(s) ds
	\le \frac{4}{(\beta-1)^2} \int_0^\infty s^{2-\beta} \psi_s^2(s) ds,
	\qquad \psi\in C_0^\infty((0,\infty)), \ \beta>1,
  \eas
  to $\beta:=\frac{2}{n}+\gamma$ and $\psi(s):=\ze(s)\vp^\frac{p}{2}(s,t)$ for $s>0, t>0, \eps\in (0,1)$,
  admissible since $\ze\in C_0^\infty((0,\infty))$ for all $\eps\in (0,1)$, shows that
  \bas
	\int_0^\infty s^{-\frac{2}{n}-\gamma} \ze^2 \vp^p ds
	&\le & \frac{4}{(\gamma-1+\frac{2}{n})^2} \int_0^\infty s^{2-\frac{2}{n}-\gamma} \cdot 
		\Big\{ \frac{p}{2} \ze \vp^\frac{p-2}{2} \vp_s + \zes \vp^\frac{p}{2}\Big\}^2 ds \nn\\
	&=&\frac{1}{(\gamma-1+\frac{2}{n})^2}\left\{p^2 \int_0^\infty s^{2-\frac{2}{n}-\gamma} \ze^2 \vp^{p-2} \vp_s^2 ds \right.\\
	& & \hspace{2.5cm}+ 4p\int_0^\infty s^{2-\frac{2}{n}-\gamma} \ze\zes \vp^{p-1} \vp_s ds \\
	& & \hspace{2.5cm}\left.\quad+ 4\int_0^\infty s^{2-\frac{2}{n}-\gamma} \zes^2 \vp^p ds\right\}
	\quad \mbox{for all $t>0$ and } \eps\in (0,1).
  \eas
  As an integration by parts shows that
  \bas
	& & \hs{-20mm}
	p\int_0^\infty s^{2-\frac{2}{n}-\gamma} \ze\zes \vp^{p-1} \vp_s ds
	+ \int_0^\infty s^{2-\frac{2}{n}-\gamma} \zes^2 \vp^p ds \\
	&=& - (2-\frac{2}{n}-\gamma)\int_0^\infty s^{1-\frac{2}{n}-\gamma} \ze\zes \vp^p ds \nn\\
	& & - \int_0^\infty s^{2-\frac{2}{n}-\gamma} (\ze\zess+\zes^2) \vp^p ds \nn\\
	& & + \int_0^\infty s^{2-\frac{2}{n}-\gamma} \zes^2 \vp^p ds \nn\\
	&=& (\gamma-2+\frac{2}{n})\int_0^\infty s^{1-\frac{2}{n}-\gamma} \ze\zes \vp^p ds \nn\\
	& & - \int_0^\infty s^{2-\frac{2}{n}-\gamma} \ze\zess \vp^p ds
	\qquad \mbox{for all $t>0$ and } \eps\in (0,1),
  \eas
  this implies \eqref{41.1}.
\qed
Let us next make sure that those integrals in (\ref{4.1}) and (\ref{41.1}) which contain $\zes$ and $\zess$ indeed vanish
in the limit $\eps\searrow 0$.
Our first statement in this regard is concerned with the corresponding cut-off process near the origin.
\begin{lem}\label{lem42}
  Let $n\ge 3$, assume (\ref{init}), and let $p>1$ and $\gamma>1-\frac{2}{n}$ be such that
  \be{42.1}
	p>\frac{n}{n-2} \gamma-1.
  \ee
  Then for each $T>0$, 
  \be{42.2}
	\int_0^T \int_{\frac{\eps}{2}}^\eps s^{1-\frac{2}{n}-\gamma} \ze|\zes| \vp^p dsdt \to 0
  \ee
  and
  \be{42.22}
	\int_0^T \int_{\frac{\eps}{2}}^\eps s^{2-\frac{2}{n}-\gamma} \zes^2 \vp^p dsdt \to 0
  \ee
  as well as
  \be{42.3}
	\int_0^T \int_{\frac{\eps}{2}}^\eps s^{2-\frac{2}{n}-\gamma} \ze|\zess| \vp^p dsdt \to 0
  \ee
  and
  \be{42.4}
	\int_0^T \int_{\frac{\eps}{2}}^\eps s^{-\gamma} \ze|\zes| \vp^{p+1} dsdt \to 0
  \ee
  as $\eps\searrow 0$.
\end{lem}
\proof
  Due to the fact that $\vp(s,t) \le 2s^{1-\frac{2}{n}}$ for all $s>0$ and $t>0$ by \eqref{pb}, a combination of
  \eqref{z2} with \eqref{z3} shows that
  \bas
	\int_0^T \int_{\frac{\eps}{2}}^\eps s^{1-\frac{2}{n}-\gamma} \ze|\zes| \vp^p dsdt
	&\le& \frac{2^p \K T}{\eps} \cdot \int_{\frac{\eps}{2}}^\eps s^{1-\frac{2}{n}-\gamma + (1-\frac{2}{n})p} ds \\
	&\le& \frac{2^p \K T}{\eps} \cdot \frac{\eps}{2} \cdot \eps^{1-\frac{2}{n}-\gamma + (1-\frac{2}{n})p} 
	\qquad \mbox{for all } \eps\in (0,1),
  \eas
  so that \eqref{42.2} becomes a consequence of \eqref{42.1}, which namely implies that
  \be{42.5}
	1-\frac{2}{n}-\gamma+\Big(1-\frac{2}{n}\Big) p
	> 1-\frac{2}{n}-\gamma + \Big(1-\frac{2}{n}\Big)\cdot \Big(\frac{n}{n-2}\gamma-1\Big)=0.
  \ee
  Similarly, 
{\cred 
  we have
  \bas
	\int_0^T \int_{\frac{\eps}{2}}^\eps s^{2-\frac{2}{n}-\gamma} \zes^2 \vp^p dsdt 
	&\le& \frac{2^p \K^2 T}{\eps^2} \cdot \frac{\eps}{2} \cdot \eps^{2-\frac{2}{n}-\gamma+(1-\frac{2}{n})p}
	\qquad \mbox{and} \\
	\int_0^T \int_{\frac{\eps}{2}}^\eps s^{2-\frac{2}{n}-\gamma} \ze|\zess| \vp^p dsdt 
	&\le& \frac{2^p \K T}{\eps^2} \cdot \frac{\eps}{2} \cdot \eps^{2-\frac{2}{n}-\gamma+(1-\frac{2}{n})p}
	\qquad \mbox{as well as} \\
	\int_0^T \int_{\frac{\eps}{2}}^\eps s^{-\gamma} \ze|\zes| \vp^{p+1} dsdt
	&\le& \frac{2^{p+1} \K T}{\eps} \cdot \frac{\eps}{2} \cdot \eps^{-\gamma+(1-\frac{2}{n})(p+1)}
  \eas
}
  for all $\eps\in (0,1)$, whence \eqref{42.22}, \eqref{42.3} and \eqref{42.4} follow upon observing that
  $-1+[2-\frac{2}{n}-\gamma+(1-\frac{2}{n})p]=1-\frac{2}{n}-\gamma+(1-\frac{2}{n})p$ and also
  $-\gamma+(1-\frac{2}{n})(p+1)=1-\frac{2}{n}-\gamma+(1-\frac{2}{n})p$,
  and again using \eqref{42.5}.
\qed
With respect to the tail regularization provided by the $\ze$,
our corresponding estimation will additionally rely on a rough local-in-time uper estimate for $\vp$ which results
from a comparison-type argument underlying the following.
\begin{lem}\label{lem3}
  Let $n\ge 3$, and suppose that besides fulfilling (\ref{init}), $u_0$ is such that with $\vp$ as in (\ref{phi}) we have
  \be{3.1}
	\sup_{s\ge 1} \frac{\vp(s,0)}{s^\al} < \infty
  \ee
  with some $\al\in (0,1-\frac{2}{n})$.
  Then there exist $\lam>0$ and $C>0$ such that
  \be{3.2}
	\vp(s,t) \le C e^{\lam t} s^\al
	\qquad \mbox{for all $s>0$ and }  t\ge 0.
  \ee
\end{lem}
\proof
  According to \eqref{3.1}, we let $c_1>0$ be such that
  \be{3.3}
	\vp(s,0) \le c_1 s^\al
	\qquad \mbox{for all } s\ge 1,
  \ee
  and set
  \be{3.4}
	\lam:=2n\al + 2(n-2).
  \ee
  Then for
  \bas
	\ophi(s,t):=y(t) s^\al,
	\qquad s\ge 1, \ t\ge 0,
  \eas
  with
  \be{3.44}
	y(t):=c_2 e^{\lam t},
	\qquad t\ge 0,
  \ee
  and with $c_2:=\max\{c_1,2\}$, from \eqref{3.3} we know that
  \be{3.5}
	\ophi(s,0) =  c_1 s^\al \ge \vp(s,0)
	\qquad \mbox{for all } s\ge 1,
  \ee
  while \eqref{3.44} together with the inequality $c_2\ge 2$ shows that due to \eqref{pb},
  \be{3.55}
	\ophi(1,t)
	= y(t) \ge c_2 \ge 2 = w_\star(1) \ge \vp(1,t)
	\qquad \mbox{for all } t>0.
  \ee
  Moreover, observing that $\ophi_s(s,t)=\al y(t) s^{\al-1} \ge 0$ and $\ophi_{ss}(s,t) = -\al(1-\al) y(t) s^{\al-2} \le 0$
  for all $s>1$ and $t>0$, thanks to the identity $y'=\lam y$ we obtain that
  \bea{3.6}
	\ophi_t - n^2 s^{2-\frac{2}{n}} \ophi_{ss} - 2n s^{1-\frac{2}{n}} \ophi_s - 2(n-2) s^{-\frac{2}{n}} \ophi + n\ophi \ophi_s
	&\ge& 	\ophi_t - 2n s^{1-\frac{2}{n}} \ophi_s - 2(n-2) s^{-\frac{2}{n}} \ophi \nn\\
	&=& y'(t) s^\al - 2n\al y(t) s^{\al-\frac{2}{n}} - 2(n-2) y(t) s^{\al-\frac{2}{n}} \nn\\
	&=& \big\{ \lam s^\frac{2}{n} - 2n\al - 2(n-2)\big\} \cdot s^{\al-\frac{2}{n}} y(t) \nn\\
	&\ge& \big\{ \lam - 2n\al - 2(n-2)\big\} \cdot s^{\al-\frac{2}{n}} y(t) \nn\\[2mm]
	&=& 0
	\qquad \mbox{for all $s>1$ and } t>0
  \eea
  because of \eqref{3.4}.
  Hence, for each $\eps\in (0,1)$ the function $d_\eps$ defined by
  \bas
	d_\eps(s,t):=\vp(s,t)-\ophi(s,t) - \eps e^{4nt} s,
	\qquad s\ge 1, \ t\ge 0,
  \eas
  satisfies
  \be{3.7}
	d_\eps(s,0) \le - \eps s < 0
	\qquad \mbox{for all } s\ge 1
  \ee
  and
  \be{3.8}
	d_\eps(1,t) \le -\eps e^{4nt} < 0
	\qquad \mbox{for all } t\ge 0
  \ee
  by \eqref{3.5} and \eqref{3.55}, as well as
  \bea{3.9}
	d_{\eps t} 
	&=&\vp_t - \ophi_t 
	- 4n \eps e^{4nt} s \nn\\
	&\le& n^2 s^{2-\frac{2}{n}} \vp_{ss} + 2n s^{1-\frac{2}{n}} \vp_s + 2(n-2) s^{-\frac{2}{n}} \vp - n\vp \vp_s \nn\\
	& & - n^2 s^{2-\frac{2}{n}} \ophi_{ss} - 2n s^{1-\frac{2}{n}} \ophi_s - 2(n-2) s^{-\frac{2}{n}} \ophi + n \ophi \ophi_s \nn\\
	& & - 4n \eps e^{4nt} s \nn\\
	&=& n^2 s^{2-\frac{2}{n}} d_{\eps ss} 
	+ 2n s^{1-\frac{2}{n}} (d_{\eps s} + \eps e^{4nt})
	+ 2(n-2) s^{-\frac{2}{n}} (d_\eps + \eps e^{4nt} s) \nn\\
	& & - n(d_\eps + \ophi+\eps e^{4nt} s)\cdot(d_{\eps s} + \ophi_s + \eps e^{4nt})
	+ n\ophi \ophi_s \nn\\
	& & - 4n \eps e^{4nt} s
	\qquad \mbox{for all $s>1$ and $t>0$}
  \eea
  according to \eqref{0p} and \eqref{3.6}.
  But since \eqref{pb} ensures that additionally
  \bas
	d_\eps(s,t) \le w_\star(s) - \eps e^{4nt} s
	\le 2s^{1-\frac{2}{n}} - \eps s < 0
	\qquad \mbox{for all $s>(\frac{2}{\eps})^\frac{n}{2}$ and } t>0,
  \eas
  from \eqref{3.7} and \eqref{3.8} it follows that for each $\eps\in (0,1)$,
  \bas
	t_\eps:=\sup \Big\{ t_0>0 \ \Big| \ d_\eps(s,t) <0 \mbox{ for all $s\ge 1$ and } t\in [0,t_0] \Big\}
  \eas
  is well-defined and positive, and that if $t_\eps$ was finite, then there would exist 
  $s_\eps\in \big(1,(\frac{2}{\eps})^\frac{n}{2}\big)$ such that $d_\eps(s_\eps,t_\eps)=0, d_{\eps t}(s_\eps,t_\eps)\ge 0$,
  $d_{\eps s}(s_\eps,t_\eps)=0$ and $d_{\eps ss}(s_\eps,t_\eps)\le 0$.
  In light of \eqref{3.9} and the nonnegativity of $\ophi$ and
  $\ophi_s$, however, this would mean that at $(s,t)=(s_\eps,t_\eps)$ we would have
  \bas
	0 
	&\le& 
	2n s^{1-\frac{2}{n}} \cdot \eps e^{4nt}
	+ 2(n-2) s^{-\frac{2}{n}} \cdot \eps e^{4nt} s \nn\\
	& & - n(\ophi+\eps e^{4nt} s) \cdot( \ophi_s + \eps e^{4nt})
	+ n\ophi \ophi_s
	- 4n \eps e^{4nt} s \\
	&\le& \big\{ 2n s^{-\frac{2}{n}} + 2(n-2)s^{-\frac{2}{n}} - 4n \big\} \cdot \eps e^{4nt} s \\
	&\le & \big\{ 2n + 2(n-2) - 4n \big\} \cdot \eps e^{4nt} s \\
	&<& 0.
  \eas
  This absurd conclusion shows that actually $d_\eps<0$ in all of $[1,\infty)\times [0,\infty)$, and that thus
  \bas
	\vp(s,t) \le c_2 e^{\lam t} s^\al
	\qquad \mbox{for all $s\ge 1$ and } t\ge 0,
  \eas
  because $\eps>0$ was arbitrary.
  As trivially
  \bas
	\frac{\vp(s,t)}{e^{\lam t} s^\al} 
	\le \frac{w_\star(s)}{e^{\lam t} s^\al}
	= 2e^{-\lam t} s^{1-\frac{2}{n}-\al} 
	\le 2
	\qquad \mbox{for all $s\in (0,1)$ and } t\ge 0
  \eas
  due to (\ref{pb}) and the fact that $\al\le 1-\frac{2}{n}$, this establishes \eqref{3.2} with
  $C:=\max\{c_2,2\}=c_2$.
\qed
If $p$ is not too large relative to $\gamma$, the latter can be used to verify that, indeed, also the approximation
errors ancountered in the considered outer parts asymptotically disappear:
\begin{lem}\label{lem43}
  Let $n\ge 3$, suppose that (\ref{init}) and \eqref{3.1} hold 
  with some $\al\in (0,1-\frac{2}{n})$, and let $p>1$ and $\gamma>1-\frac{2}{n}$ satisfy
  \be{43.1}
	p\al<\gamma-1+\frac{2}{n}.
  \ee
  For arbitrary $T>0$, as $\eps\searrow 0$ it then follows that
  \be{43.2}
	\int_0^T \int_{\frac{1}{\eps}}^{\frac{2}{\eps}} s^{1-\frac{2}{n}-\gamma} |\ze\zes| \vp^p dsdt \to 0,
  \ee
  that
  \be{43.22}
	\int_0^T \int_{\frac{1}{\eps}}^{\frac{2}{\eps}} s^{2-\frac{2}{n}-\gamma} \zes^2 \vp^p dsdt \to 0,
  \ee
  and that
  \be{43.3}
	\int_0^T \int_{\frac{1}{\eps}}^{\frac{2}{\eps}} s^{2-\frac{2}{n}-\gamma} |\ze\zess| \vp^p dsdt \to 0.
  \ee
%
\end{lem}
\proof
  In view of Lemma \ref{lem3}, we can fix $c_1=c_1(T)>0$ such that
  \bas
	\vp(s,t) \le c_1 s^\al
	\qquad \mbox{for all $s>0$ and } t\in (0,T),
  \eas
  whence applications of \eqref{z2} and \eqref{z4} yield the inequalities
{\cred
  \bas
	\int_0^T \int_{\frac{1}{\eps}}^{\frac{2}{\eps}} s^{1-\frac{2}{n}-\gamma} |\ze\zes| \vp^p dsdt 
	&\le&	2^{1-\frac{2}{n}-\gamma+p\al} c_1^p \K T \cdot \eps^{-1+\frac{2}{n}+\gamma-p\al}
	\qquad \mbox{and} \\
	\int_0^T \int_{\frac{1}{\eps}}^{\frac{2}{\eps}} s^{2-\frac{2}{n}-\gamma} \zes^2 \vp^p dsdt
	&\le&	2^{2-\frac{2}{n}-\gamma+p\al} c_1^p \K^2 T \cdot \eps^{-1+\frac{2}{n}+\gamma-p\al}
	\qquad \mbox{as well as} \\
	\int_0^T \int_{\frac{1}{\eps}}^{\frac{2}{\eps}} s^{2-\frac{2}{n}-\gamma} |\ze\zess| \vp^p dsdt 
	&\le&	2^{2-\frac{2}{n}-\gamma+p\al} c_1^p \K T \cdot \eps^{-1+\frac{2}{n}+\gamma-p\al}
  \eas
}
  for all $\eps\in (0,1)$. Since from \eqref{43.1} we know that
  $-1+\frac{2}{n}+\gamma-p\al>0$, \eqref{43.2}, \eqref{43.22}, and \eqref{43.3} follow.
\qed
For choices of $p$ compatible with the requirements in Lemma \ref{lem42} and Lemma \ref{lem43},
we can now accomplish our main step toward Theorem \ref{theo7} on the basis of \eqref{41.1}
and the circumstance that the third summand on the right of \eqref{4.1} is favorably signed.
\begin{lem}\label{lem44}
  Let $n\ge 10$, $\al\in (0,1-\frac{2}{n})$ and $\gamma>1-\frac{2}{n}$, and let $p>1$ be such that
  \eqref{42.1} holds, that
  \be{44.01}
	p\al<\gamma-1,
  \ee
  and that
  \be{44.1}
	p_-(\gamma) \le p \le p_+(\gamma),
  \ee
  where
  \be{44.2}
	p_\pm(\gamma):=\frac{n}{4}\cdot \Big(\gamma-1+\frac{2}{n}\Big) \cdot \Big(1 \pm \sqrt{\frac{n-10}{n-2}}\Big).
  \ee
  Then whenever $u_0$ satisfies (\ref{init}) and \eqref{3.1}, it follows that
  \be{44.3}
	\int_0^\infty \int_0^\infty s^{-\gamma-1} \vp^{p+1}(s,t) dsdt < \infty.
  \ee
\end{lem}
\proof
  A combination of Lemma \ref{lem4} with Lemma \ref{lem41} provides $c_1>0$ such that the functions in \eqref{pe} satisfy
  \bea{44.4}
	\frac{1}{p}\phi_\eps'(t)
	&\le& - n^2(p-1) \int_0^\infty s^{2-\frac{2}{n}-\gamma} \ze^2 \vp^{p-2} \vp_s^2 ds \nn\\
	& & + \Big\{ \frac{n}{p} \cdot (n\gamma-2n+4)\cdot \Big(\gamma-1+\frac{2}{n}\Big) + 2(n-2)\Big\}
	\cdot \frac{p^2}{(\gamma-1+\frac{2}{n})^2} \int_0^\infty s^{2-\frac{2}{n}-\gamma} \ze^2 \vp^{p-2} \vp_s^2 ds \nn\\
	& & + c_1 \int_0^\infty s^{1-\frac{2}{n}-\gamma} \ze|\zes| \vp^p ds
	+ c_1 \int_0^\infty s^{2-\frac{2}{n}-\gamma} \zes^2 \vp^p ds \nn\\
	& & + c_1 \int_0^\infty s^{2-\frac{2}{n}-\gamma} \ze |\zess| \vp^p ds \nn\\
	& & - \frac{n\gamma}{p+1} \int_0^\infty s^{-\gamma-1} \ze^2 \vp^{p+1} ds
	+ \frac{2n}{p+1} \int_0^\infty s^{-\gamma} \ze\zes \vp^{p+1} ds
  \eea
  for all $t>0$ and $\eps\in (0,1)$.
  Here we note that according to \eqref{44.2} and \eqref{44.1},
  \bas
	& & \hs{-10mm}
	- n^2(p-1) 
	+ \Big\{ \frac{n}{p} \cdot (n\gamma-2n+4)\cdot \Big(\gamma-1+\frac{2}{n}\Big) + 2(n-2)\Big\}
	\cdot \frac{p^2}{(\gamma-1+\frac{2}{n})^2} \\
	&=& \frac{2(n-2)}{(\gamma-1+\frac{2}{n})^2} \cdot \bigg\{
	- \frac{n^2(\gamma-1+\frac{2}{n})^2}{2(n-2)} \cdot (p-1)
	+ \frac{n}{2(n-2)} \cdot (n\gamma-2n+4)\cdot \Big(\gamma-1+\frac{2}{n}\Big) \cdot p
	+ p^2 \bigg\} \\
	&=& \frac{2(n-2)}{(\gamma-1+\frac{2}{n})^2} \cdot \big(p-p_-(\gamma)\big)\cdot \big(p-p_+(\gamma)\big) \le 0
  \eas
  so that an integration in \eqref{44.4} shows that for all $T>0$ and $\eps\in (0,1)$,
  \bea{44.5}
	\frac{1}{p}\phi_\eps(T) 
	+ \frac{n\gamma}{p+1} \int_0^T \int_0^\infty s^{-\gamma-1} \ze^2 \vp^{p+1} dsdt
	&\le& \frac{1}{p}\phi_\eps(0) 
	+ c_1 \int_0^T \int_0^\infty s^{1-\frac{2}{n}-\gamma} |\ze\zes| \vp^p dsdt \nn\\
	& & + c_1 \int_0^T \int_0^\infty s^{2-\frac{2}{n}-\gamma} \zes^2 \vp^p dsdt \nn\\
	& & + c_1 \int_0^T \int_0^\infty s^{2-\frac{2}{n}-\gamma} |\ze\zess| \vp^p dsdt \nn\\
	& & + \frac{2n}{p+1} \int_0^T \int_0^\infty s^{-\gamma} \ze\zes \vp^{p+1} dsdt.
  \eea
  Now as we are assuming \eqref{42.1} and \eqref{44.01}, with the latter evidently being sharper than \eqref{43.1}, we may rely on 
  Lemma \ref{lem42} and Lemma \ref{lem43} to see that since $\zes\equiv 0$ 
  outside $[\frac{\eps}{2},\eps]\cup [\frac{1}{\eps},\frac{2}{\eps}]$ for all $\eps\in (0,1)$ by (\ref{z2}),
  for each fixed $T>0$ we have
  \bas
	c_1 \int_0^T \int_0^\infty s^{1-\frac{2}{n}-\gamma} |\ze\zes| \vp^p dsdt
	+ c_1 \int_0^T \int_0^\infty s^{2-\frac{2}{n}-\gamma} \zes^2 \vp^p dsdt 
	+ c_1 \int_0^T \int_0^\infty s^{2-\frac{2}{n}-\gamma} |\ze\zess| \vp^p dsdt 
	\to 0
  \eas
  and
  \bas
	\frac{2n}{p+1} \int_0^T \int_0^1 s^{-\gamma} \ze\zes \vp^{p+1} dsdt
	\to 0
  \eas
  as $\eps\searrow 0$.
  Since \eqref{chi} and \eqref{ze} ensure that $\zes\le 0$ on $(1,\infty)$ for all $\eps\in (0,1)$, and since clearly
  $\ze\nearrow 1$ in $(0,\infty)$ as $\eps\searrow 0$ by \eqref{z2}, from \eqref{44.5} and Beppo Levi's theorem we thus infer that
  \bas
	\frac{n\gamma}{p+1} \int_0^T \int_0^\infty s^{-\gamma-1} \vp^{p+1} dsdt
	\le \liminf_{\eps\searrow 0} \frac{1}{p} \phi_\eps(0)
	= \frac{1}{p} \int_0^\infty s^{-\gamma} \vp^p(s,0) ds
	\qquad \mbox{for all } T>0.
  \eas
  It thus remains to observe that the latter integral is finite, which indeed results from (\ref{pb}), (\ref{3.2}) and the fact that
  \bas
	\int_0^1 s^{-\gamma} \cdot s^{(1-\frac{2}{n})p} ds
	+ \int_1^\infty s^{-\gamma} \cdot s^{p\al} ds 
	< \infty
  \eas
  according to the inequalities $\frac{n-2}{n}p > \gamma-1$ and $p\al<\gamma-1$ guaranteed by \eqref{42.1} and \eqref{44.01}.
\qed
According to appropriate parabolic smoothing effects, the weak decay property expressed in \eqref{44.3}
can be turned into the following statement on genuine vanishing of $\vp$ in the large time limit.
\begin{lem}\label{lem45}
  Let $n\ge 10$, $\al\in (0,1-\frac{2}{n})$, $\gamma>1-\frac{2}{n}$ and $p>1$ satisfy
  \eqref{42.1}, \eqref{44.01} and \eqref{44.1}, and 
  suppose that $u_0$ complies with (\ref{init}) and \eqref{3.1}.
  Then
	\[
	\vp(\cdot,t)\to 0
	\quad \mbox{in } C^1_{loc}((0,\infty))
	\qquad \mbox{as } t\to\infty.
	\]
\end{lem}
\proof
  If the claim was false, then there would exist $s_1>0, s_2>s_1$, $c_1>0$ and $(t_k)_{k\in\N}\subset (1,\infty)$ such that
  $t_k\to\infty$ as $k\to\infty$ and
  \be{45.2}
	\|\vp(\cdot,t_k)\|_{C^1([s_1,s_2])} \ge c_1
	\qquad \mbox{for all } k\in\N.
  \ee
  Keeping these selections of $s_1$ and $s_2$ fixed, we see that the parabolic equation in \eqref{0p} is of the form
  \be{45.3}
	\vp_t = D(s)\vp_{ss} + A(s,t)\vp_s + B(s)\vp,
	\qquad (s,t)\in (0,\infty)^2,
  \ee
  where $D(s):=n^2 s^{2-\frac{2}{n}}$, $A(s,t):=2ns^{1-\frac{2}{n}} - n\vp(s,t)$ and $B(s):=2(n-2)s^{-\frac{2}{n}}$,
  $s\in (0,\infty), \ t>0$, satisfy
  $0<n^2 \cdot (\frac{s_1}{2})^{2-\frac{2}{n}} \le D(s) \le n^2 \cdot (2s_2)^{2-\frac{2}{n}}$ and, by (\ref{pb}),
  $0\le A(s,t) \le 2n \cdot (2s_2)^{1-\frac{2}{n}}$ as well as $0\le B(s) \le 2(n-2) \cdot (\frac{s_1}{2})^{-\frac{2}{n}}$
  for all $s\in [\frac{s_1}{2},2s_2]$ and $t>0$.
  Since $\vp$ is bounded in $[\frac{s_1}{2},2s_2] \times (0,\infty)$ by (\ref{pb}), we may therefore apply interior parabolic 
  Schauder theory (\cite{LSU}) to the non-degenerate equation \eqref{45.3} to infer the existence of $\vt\in (0,1)$ and $c_2>0$
  such that
  \be{45.4}
	\|\vp\|_{C^{2+\vt,1-\frac{\vt}{2}}([s_1,s_2]\times [t,t+1])} \le c_2
	\qquad \mbox{for all } t\ge 1,
  \ee
  which thanks to the Arzel\`a-Ascoli theorem in particular implies that $(\vp(\cdot,t_k))_{k\in\N}$ is relatively compact in
  $C^1([s_1,s_2])$.
  We can thus find $\vp_\infty\in C^1([s_1,s_2])$ and a subsequence $(t_{k_l})_{l\in\N}$ of $(t_k)_{k\in\N}$ such that
  \be{45.5}
	\vp(\cdot,t_{k_l}) \to \vp_\infty
	\quad \mbox{in } C^1([s_1,s_2])
	\qquad \mbox{as } l\to\infty,
  \ee
  where \eqref{45.2} asserts that $\vp_\infty\not\equiv 0$, and that thus there exist $s_3\in [s_1,s_2)$, $s_4\in (s_3,s_2]$ and
  $c_3>0$ fulfilling
  $\vp_\infty \ge c_3$ in $[s_3,s_4]$.
  In line with \eqref{45.5}, we can therefore choose $l_0\in\N$ such that $\vp(\cdot,t_{k_l}) \ge \frac{c_3}{2}$
  in $[s_3,s_4]$ for all $l\ge l_0$, whence observing that \eqref{45.4} moreover warrants uniform continuity of $\vp$ in
  $[s_1,s_2]\times [1,\infty)$, we can pick $\del>0$ such that, in fact, 
  \bas
	\vp(s,t)\ge \frac{c_3}{4}
	\qquad \mbox{for all $s\in [s_3,s_4]$, $t\in [t_{k_l},t_{k_l}+\del]$ and } l\ge l_0.
  \eas
  This, however, means that
  \bas
	\int_0^\infty \int_0^\infty s^{-\gamma-1} \vp^{p+1} (s,t) dsdt
	\ge \Big(\frac{c_3}{4}\Big)^{p+1} \cdot s_4^{-\gamma-1} \cdot 
		\Bigg| \bigcup_{l_0 \le l \le N} [t_{k_l},t_{k_l}+\del] \Bigg|
	\qquad \mbox{for all } N\ge l_0
  \eas
  and thereby contradicts \eqref{44.3}, because here the expression on the right-hand side diverges to $+\infty$ as $N\to\infty$.
\qed
It remains to suitably adjust parameters in order to establish our main result on stability and attractiveness of the 
singular equilibrium from (\ref{stat}) in the high-dimensional setting under consideration:\abs
\proofc of Theorem \ref{theo7}. \quad
  Since $\frac{n-2+\sqrt{(n-2)(n-10)}}{2}<\frac{n-2+\sqrt{(n-2)(n-2)}}{2}=n-2$, we may assume without loss of generality that
  $\theta<n-2$.
  Then $\al:=\frac{n-2-\theta}{n}$ is positive and clearly satisfies $\al<1-\frac{2}{n}$, and thanks to (\ref{theta}) it follows
  that moreover
  \be{7.5}
	\al < \frac{n-2-\frac{n-2+\sqrt{(n-2)(n-10)}}{2}}{n}
	= \frac{n-2-\sqrt{(n-2)(n-10)}}{2n}.
  \ee
  On the other hand, it can easily be verified that
  \bas
	\frac{\gamma-1}{\frac{n}{4}(\gamma-1+\frac{2}{n})(1+\sqrt{\frac{n-10}{n-2}})}
	&\to& \frac{1}{\frac{n}{4}(1+\sqrt{\frac{n-10}{n-2}})}
	=\frac{n-2-\sqrt{(n-2)(n-10)}}{2n}
	\qquad \mbox{as } \gamma\to\infty,
  \eas
  so that, since furthermore the 
{\cred 
  assumption $n\ge 10$ warrants that
}
  \bas
	\frac{\frac{n}{n-2}\gamma -1}{\frac{n}{4}(\gamma-1+\frac{2}{n})(1+\sqrt{\frac{n-10}{n-2}})}
	\to \frac{\frac{n}{n-2}}{\frac{n}{4}(1+\sqrt{\frac{n-10}{n-2}})}
	 < 1
	\qquad \mbox{as } \gamma\to\infty,
  \eas
  due to \eqref{7.5} we can choose $\gamma>1-\frac{2}{n}$ suitably large such that
  \be{7.6}
	\frac{\gamma-1}{\frac{n}{4}(\gamma-1+\frac{2}{n})(1+\sqrt{\frac{n-10}{n-2}})}
	>\al,
  \ee
  and that
  \be{7.7}
	\frac{\frac{n}{n-2}\gamma -1}{\frac{n}{4}(\gamma-1+\frac{2}{n})(1+\sqrt{\frac{n-10}{n-2}})}
	< 1,
  \ee
  where enlarging $\gamma$ if necessary we can also achieve that
  \be{7.8}
	\frac{n}{4}\cdot\Big(\gamma-1+\frac{2}{n}\Big)\cdot\Big(1+\sqrt{\frac{n-10}{n-2}}\Big) >1.
  \ee
  If now we define
  \be{7.9}
	p:=\frac{n}{4}\cdot\Big(\gamma-1+\frac{2}{n}\Big)\cdot\Big(1+\sqrt{\frac{n-10}{n-2}}\Big),
  \ee
  then from \eqref{7.8} we obtain that $p>1$, while \eqref{7.7} and \eqref{7.6} assert that $p>\frac{n}{n-2}\gamma-1$
  and $p\al<\gamma-1$.\abs
  To complete our preparation for an application of Lemma \ref{lem45}, we now note that according the inequality
  $\theta<n-2$ and \eqref{7.2}, the function $\vp$ from \eqref{phi} has the property that
  \bas
	\vp(s,0)
	&=& w_\star(s) - \int_0^{s^\frac{1}{n}} \rho^{n-1} u_0(\rho) d\rho \\
	&=& 2s^{1-\frac{2}{n}}
	- \int_0^1 \rho^{n-1} u_0(\rho)d\rho
	- \int_1^{s^\frac{1}{n}} \rho^{n-1} \cdot \Big\{ \frac{2(n-2)}{\rho^2} - \frac{C}{\rho^{2+\theta}} \Big\} d\rho \\
	&\le& 2 + \frac{C}{n-2-\theta} \cdot s^\frac{n-2-\theta}{n} \\
	&<& \Big(2+\frac{C}{n-2-\theta}\Big) \cdot s^\al
	\qquad \mbox{for all } s>1.
  \eas
  Since thus \eqref{3.1} is satisfied, and since (\ref{7.9}) together with \eqref{44.2} ensures that also \eqref{44.1} holds,
  we may therefore indeed employ Lemma \ref{lem45} to conclude that $\vp(\cdot,t)\to 0$ in $C^1_{loc}((0,\infty))$ as $t\to\infty$.
  Due to the identity $u(r,t)=n(w_\star)_s(r^n,t) - n\vp_s(r^n,t)$ for $r>0$ and $t>0$, this is equivalent to
  \eqref{7.4}, whereupon \eqref{7.3} becomes an evident consequence of the fact that $u_\star$ is unbounded.
\qed
\mysection{Instability of $u_\star$. Proofs of Theorem \ref{theo25} and Proposition \ref{prop26}}\label{sect4}
{\cred
In contrast to the above analysis, our considerations related to the announced instability results
will be based on comparison arguments, at their core relying on favorable parabolic inequalities satisfied by some
essentially explicit functions of separated structure that quantify deviations from equilibrium.\abs
The core of our approach in this regard can be found in the following lemma
which in low-dimensional spaces identifies a fundamental repulsion property of $w_\star$, here yet observed at certain particular
points in space.
}
The proof develops a basic qualitative strict positivity feature of $\vp$, as directly resulting from a strong maximum principle,
into a spatially local but quantitative lower estimate for this deviation by means of a comparison argument drawing
on a corresponding instability feature of the trivial solution to a Bernoulli-type ODE.
\begin{lem}\label{lem21}
  Let $3\le n\le 9$.
  Then for all $s_0>0$ there exist $s_\star=s_\star(s_0)>s_0$ and
  $M=M(s_0)>0$ with the property that whenever $u_0$ satisfies (\ref{init}), one can find $t_0=t_0(u_0)>0$ such that
  for the function $w$ from \eqref{w} we have
	\[
	w(s_\star,t)\le w_\star(s_\star) - M
	\qquad \mbox{for all } t\ge t_0.
	\]
\end{lem}
\proof
  Using that $3\le n\le 9$, we fix $\om>0$ such that $\om^2 < \frac{(n-2)(10-n)}{4n^2}$, and note that if we furthermore abbreviate
  $\mu:=\frac{n-2}{2n}$, then this implies that
  \be{21.2}
	c_1:=n^2(\mu^2-\mu-\om^2) + 2n\mu + 2(n-2)
  \ee
  satisfies
  \bas
	c_1
	&>& n^2 \cdot \Big\{ \mu^2-\mu-\frac{(n-2)(10-n)}{4n^2}\Big\} + 2n\mu + 2(n-2) 	\\
	&=& \frac{n-2}{4} \cdot \big\{ n-2-2n-(10-n)+12\big\} = 0.
  \eas
  With this choice fixed, given $s_0>0$ we let $k_0=k_0(s_0)\in\N$ be large enough such that $e^\frac{(4k_0-1)\pi}{2\om}>s_0$,
  and define
  \bas
	s_1\equiv s_1(s_0):=e^\frac{(4k_0-1)\pi}{2\om},
	\quad 
	s_2\equiv s_2(s_0):=e^\frac{(4k_0+1)\pi}{2\om}
	\quad \mbox{and} \quad
	s_\star\equiv s_\star(s_0):=e^\frac{2k_0\pi}{\om},
  \eas
  so that $s_0<s_1<s_\star<s_2$, with
  \be{21.3}
	\cos (\om\ln s)
	\lball
	>0 \qquad  \mbox{for all } s\in (s_1,s_2), \\[1mm]
	=0 \qquad \mbox{if } s\in \{s_1,s_2\}, \\[1mm]
	=1 \qquad \mbox{if } s=s_\star.
	\ear
  \ee
  We moreover set
  \be{21.4}
	c_2\equiv c_2(s_0):=s_2^\frac{2}{n}
	\qquad \mbox{and} \qquad
	c_3\equiv c_3(s_0):=n(\mu+\om) s_1^{-\mu},
  \ee
  and apply the classical parabolic strong maximum principle (\cite{hulshof}) to (\ref{0p}) to see that $\vp>0$
  in $[s_1,s_2]\times (0,\infty)$, whence in particular we can find $c_4=c_4(u_0)>0$ fulfilling
  \be{21.5}
	\vp(s,1)\ge c_4
	\qquad \mbox{for all } s\in [s_1,s_2].
  \ee
  The positivity of all the $c_i$ for $i\in\{1,2,3,4\}$ now ensures that the solution $y\in C^1([1,\infty))$ of the Bernoulli type
  problem
  \be{21.6}
	\lbal
	y'(t)= \frac{c_1}{c_2} y(t) - \frac{c_3}{c_2} y^2(t),
	\qquad t>1, \\[1mm]
	y(1)=c_4 s_2^{-\mu},
	\ear
  \ee
  satisfies $y(t)\to \frac{c_1}{c_3}>0$ as $t\to\infty$, so that there exists $t_0=t_0(u_0)>1$ such that
  \bas
	y(t) \ge \frac{c_1}{2c_3}
	\qquad \mbox{for all } t>t_0.
  \eas
  The function $\uphi\in C^\infty([s_1,s_2]\times [1,\infty))$ defined by letting
  \be{21.7}
	\uphi(s,t):=y(t) s^\mu \cos (\om\ln s),
	\qquad s\in [s_1,s_2], \ t\ge 1,
  \ee
  thus has the property that due to (\ref{21.3}),
  \bas
	\uphi(s_\star,t) = y(t) s_\star^\mu
	\ge \frac{c_1 s_\star^\mu}{2c_3}
	\qquad \mbox{for all } t>t_0,
  \eas
  so that since the rightmost expression here depends on $s_0$ only but not on $u_0$, the claim will follow as soon as we have
  shown that
  \be{21.8}
	\vp(s,t) \ge \uphi(s,t)
	\qquad \mbox{for all $s\in [s_1,s_2]$ and } t\ge 1.
  \ee
  To verify this, we observe that by (\ref{21.7}), (\ref{21.3}) and (\ref{pb}),
  \be{21.9}
	\uphi(s_i,t)=0 \le \vp(s_i,t)
	\qquad \mbox{for all $t\ge 1$ and for } i\in\{1,2\},
  \ee
  while due to (\ref{21.7}), (\ref{21.6}) and (\ref{21.5}),
  \be{21.10}
	\uphi(s,1) \le y(1) s^\mu \le y(1) s_2^\mu = c_4 \le \vp(s,1)
	\qquad \mbox{for all } s\in [s_1,s_2].
  \ee
  Moreover, computing
  \bas
	\uphi_t(s,t)=y'(t) s^\mu \cos(\om\ln s)
	\qquad \mbox{and} \qquad
	\uphi_s(s,t)=\mu y(t) s^{\mu-1} \cos(\om\ln s) - \om y(t) s^{\mu-1} \sin(\om\ln s)
  \eas
  as well as
  \bas
	\uphi_{ss}(s,t)
	&=& \mu(\mu-1) y(t) s^{\mu-2} \cos(\om\ln s) - \mu\om y(t) s^{\mu-2} \sin(\om\ln s) \\
	& & - (\mu-1) \om y(t) s^{\mu-2} \sin(\om\ln s) - \om^2 y(t) s^{\mu-2} \cos(\om\ln s) \\
	&=& (\mu^2-\mu-\om^2) y(t) s^{\mu-2} \cos(\om\ln s) - (2\mu-1) \om y(t) s^{\mu-2} \sin(\om\ln s)
  \eas
  for $s\in (s_1,s_2)$ and $t>1$, we see that
  \bea{21.11}
	& & \hs{-20mm}
	\uphi_t - n^2 s^{2-\frac{2}{n}} \uphi_{ss} - 2ns^{1-\frac{2}{n}} \uphi_s - 2(n-2) s^{-\frac{2}{n}} \uphi + n\uphi \uphi_s 
		\nn\\
	&=& y'(t) s^\mu \cos(\om\ln s) \nn\\
	& & - n^2 s^{2-\frac{2}{n}} \cdot \Big\{ (\mu^2-\mu-\om^2) y(t) s^{\mu-2} \cos(\om\ln s)
	- (2\mu-1)\om y(t) s^{\mu-2} \sin(\om\ln s) \Big\} \nn\\
	& & - 2n s^{1-\frac{2}{n}} \cdot \Big\{ \mu y(t) s^{\mu-1} \cos(\om\ln s)
	- \om y(t) s^{\mu-1} \sin(\om\ln s) \Big\} \nn\\
	& & -2(n-2) s^{-\frac{2}{n}} \cdot y(t) s^\mu \cos(\om\ln s) \nn\\
	& & + ny(t) s^\mu \cos(\om\ln s) \cdot 		
	\Big\{ \mu y(t) s^{\mu-1} \cos(\om\ln s) - \om y(t) s^{\mu-1} \sin(\om\ln s) \Big\} \nn\\
	&=& y(t) s^{\mu-\frac{2}{n}} \cos(\om\ln s) \cdot \Big\{
	\frac{y'(t)}{y(t)} \cdot s^\frac{2}{n} \nn\\
	& & \hs{10mm}
	- n^2(\mu^2-\mu-\om^2) - 2n\mu - 2(n-2) \nn\\
	& & \hs{10mm}
	+ n y(t) s^{\mu-1+\frac{2}{n}} \cdot \big[ \mu \cos(\om\ln s) - \om \sin(\om\ln s) \big] \Big\} \nn\\
	& & + y(t) s^{\mu-\frac{2}{n}} \sin(\om\ln s) \cdot \Big\{ n^2(2\mu-1) \om + 2n\om \Big\} \nn\\
	&=& y(t) s^{\mu-\frac{2}{n}} \cos(\om\ln s) \cdot \Big\{
	\frac{y'(t)}{y(t)} \cdot s^\frac{2}{n}
	- c_1 + ny(t) \cdot s^{-\mu} \cdot \big[ \mu \cos (\om\ln s) - \om \sin(\om\ln s) \big] \Big\}
  \eea
  for all $s\in (s_1,s_2)$ and $t>1$ by (\ref{21.2}), because our choice of $\mu$ ensures that
  \bas
	\mu-1+\frac{2}{n}=\frac{n-2}{2n}-1+\frac{2}{n} = -\frac{1}{2}+\frac{1}{n}=-\mu,
  \eas
  and that
  \bas
	n^2(2\mu-1)\om + 2n\om = n^2\cdot \Big(\frac{n-2}{n}-1\Big)\om  + 2n\om =0.
  \eas
  As using (\ref{21.4}) and (\ref{21.6}) we see that
  \bas
	\frac{y'(t)}{y(t)} \cdot s^\frac{2}{n} - c_1 + ny(t)\cdot s^{-\mu} \cdot \big[ \mu\cos(\om\ln s) -\om \sin(\om\ln s)\big]
	&\le& c_2 \cdot \frac{y'(t)}{y(t)}
	- c_1 + c_3 y(t) \nn\\[2mm]
	&=& 0
	\qquad \mbox{for all } t>1,
  \eas
  from (\ref{21.11}) we infer that $\uphi$ is a subsolution of (\ref{0p}), whence in view of (\ref{21.9}) and (\ref{21.10})
  we may rely on classical parabolic comparison (\cite{hulshof}) to conclude that indeed (\ref{21.8}) holds,
  and that thus the proof can be completed as intended.
\qed
By a second maximum-principle-based argument involving time-dependent relatives of some comparison functions already used in
\cite{biler_nadzieja1995}, the outcome of Lemma \ref{lem21} can further be developed so as to yield spatially uniform
repulsion throughout any ball of fixed radius:
\begin{lem}\label{lem22}
  Let $3\le n\le 9$.
  Then for all $s_0>0$ there exists $B=B(s_0)>0$ such that
  if (\ref{init}) holds, then with some $t_1=t_1(u_0)>0$ we have
  \be{22.1}
	w(s,t) \le \frac{2s}{s^\frac{2}{n}+B}
	\qquad \mbox{for all $s\in [0,s_0]$ and } t\ge t_1.
  \ee
\end{lem}
\proof
  Given $s_0>0$ we let $s_\star=s_\star(s_0)>s_0$ and $M=M(s_0)>0$ be as in Lemma \ref{lem21}, and take $B=B(s_0)>0$ small enough
  such that
  \be{22.2}
	4B s_\star^{1-\frac{4}{n}} \le M.
  \ee
  To see that this number has the claimed properties, we fix $u_0$ fulfilling (\ref{init}), to obtain from (\ref{w}) that
  if we take $t_0=t_0(u_0)>0$ from Lemma \ref{lem21}, and use that 
  $c_1\equiv c_1(u_0):=\|u(\cdot,t_0)\|_{L^\infty(\R^n)}$ is finite by Proposition \ref{prop0}, then
  \be{22.3}
	w(s,t_0) \le \frac{c_1 s}{n}
	\qquad \mbox{for all } s\ge 0.
  \ee
  Moreover, an application of the strong maximum principle (\cite{hulshof}) to (\ref{0p}) readily shows that if we let
  \be{22.33}
	s_1\equiv s_1(u_0):=\Big(\frac{n}{c_1}\Big)^\frac{n}{2},
  \ee
  then in the case when $s_1<s_\star$ we have $w(\cdot,t_0)<w_\star$ in $[s_1,s_\star]$, so that there exists $c_2=c_2(u_0)>0$
  fulfilling
  \be{22.4}
	w(s,t_0) \le 2 s^{1-\frac{2}{n}} - c_2
	\qquad \mbox{if } s\in [s_1,s_\star].
  \ee
  We then take $b_0=b_0(u_0)>0$ suitably small such that
  \be{22.5}
	b_0<2B,
	\qquad b_0\le \frac{n}{c_1}
	\qquad \mbox{and} \qquad
	2b_0 \max\{s_\star^{1-\frac{4}{n}},s_1^{1-\frac{4}{n}}\} \le c_2,
  \ee
  and let $b\in C^1([t_0,\infty))$ denote the solution of 
  \be{22.6}
	\lbal
	b'(t) = \frac{2b(t)}{s_\star^\frac{2}{n} + b(t)} \cdot \frac{2B-b(t)}{2B},
	\qquad t>t_0, \\[2mm]
	b(t_0)=b_0,
	\ear
  \ee
  noting that since $b_0<2B$ by (\ref{22.5}), the function $b$ increases on $(t_0,\infty)$ and satisfies
	\[
	b(t) \to 2B
	\qquad \mbox{as } t\to\infty,
	\]
  whence with some $t_1=t_1(u_0)>t_0$,
  \be{22.7}
	b(t) \ge B
	\qquad \mbox{for all } t>t_1.
  \ee
  Now, writing
	\[
	\ow(s,t) := \frac{2s}{s^\frac{2}{n} + b(t)},
	\qquad s\in [0,s_\star], \ t\ge t_0,
	\]
  we compute
  \bas
	\ow_t(s,t)= - \frac{2b'(t) s}{(s^\frac{2}{n}+b(t))^2}
  \eas
  and 
  \bas
	\ow_s(s,t)
	= \frac{2(s^\frac{2}{n}+b(t)) - 2s\cdot\frac{2}{n} s^{\frac{2}{n}-1}}{(s^\frac{2}{n}+b(t))^2}
	= \frac{(2-\frac{4}{n})s^\frac{2}{n} + 2b(t)}{(s^\frac{2}{n}+b(t))^2}
  \eas
  as well as
  \bas
	\ow_{ss}(s,t)
	&=& \frac{-\frac{4}{n}(1-\frac{2}{n}) s^{-1+\frac{4}{n}} -(\frac{4}{n}+\frac{8}{n^2}) b(t) s^{-1+\frac{2}{n}}}
		{(s^\frac{2}{n}+b(t))^3}
  \eas
  for $(s,t)\in (0,s_\star)\times (t_0,\infty)$, so that
  \bea{22.9}
	\ow_t - n^2 s^{2-\frac{2}{n}} \ow_{ss} - n\ow\ow_s
	&=& \frac{1}{(s^\frac{2}{n}+b(t))^3} \cdot \bigg\{
	-2s(s^\frac{2}{n}+b(t))b'(t) \nn\\
	& & \hs{10mm}
	+ n^2 s^{2-\frac{2}{n}} \cdot \frac{4}{n}\Big(1-\frac{2}{n}\Big) s^{-1+\frac{4}{n}}
	+ n^2 s^{2-\frac{2}{n}} \cdot \Big(\frac{4}{n}+ \frac{8}{n^2}\Big) b(t) s^{-1+\frac{2}{n}} \nn\\
	& & \hs{10mm}
	- n\cdot 2s\cdot \Big\{ \Big(2-\frac{4}{n}\Big) s^\frac{2}{n} + 2b(t)\Big\} \bigg\} \nn\\
	&=& \frac{1}{(s^\frac{2}{n}+b(t))^3} \cdot \Big\{
	-2s^{1+\frac{2}{n}} b'(t) + (8-2b'(t)) sb(t)\Big\}
  \eea
  for all $s\in (0,s_\star)$ and $t>t_0$.
  Since (\ref{22.6}) particularly entails that
  \bas
	b'(t) \le \frac{2b(t)}{b(t)} \cdot \frac{2B}{2B} =2
	\qquad \mbox{for all } t>t_0,
  \eas
  and that moreover also
  \bas
	2s^{1+\frac{2}{n}} b'(t)
	\le 2s s_\star^\frac{2}{n} \cdot \frac{2b(t)}{s_\star^\frac{2}{n}} \cdot \frac{2B}{2B}
	\le 4sb(t)
	\qquad \mbox{for all $s\in (0,s_\star)$ and } t>t_0,
  \eas
  we can here estimate
  \bas
	- 2s^{1+\frac{2}{n}} b'(t)
	+ (8-2b'(t)) sb(t)
	\ge -4sb(t) + (8-2\cdot 2) sb(t)=0
	\qquad \mbox{for all $s\in (0,s_\star)$ and } t>t_0,
  \eas
  so that from (\ref{22.9}) we obtain that
  \be{22.10}
	\ow_t - n^2 s^{2-\frac{2}{n}} \ow_{ss} - n\ow\ow_s \ge 0
	\qquad \mbox{for all $s\in (0,s_\star)$ and } t>t_0.
  \ee
  Apart from that, we clearly have
  \be{22.11}
	\ow(0,t) = 0 = w(0,t)
	\qquad \mbox{for all } t>t_0,
  \ee
  while an application of Lemma \ref{lem21} shows that thanks to (\ref{22.2}),
  \bea{22.12}
	\ow(s_\star,t)-w(s_\star,t)
	&\ge& \frac{2s_\star}{s_\star^\frac{2}{n}+b(t)}
	- \big( 2s_\star^{1-\frac{2}{n}}- M\big) 
	\ge - 2b(t) s_\star^{1-\frac{4}{n}} + M \nn\\
	&\ge& - 4B s_\star^{1-\frac{4}{n}} + M \ge 0
	\qquad \mbox{for all } t>t_0.
  \eea
  Moreover, a combination of (\ref{22.3}) with (\ref{22.5}) and (\ref{22.33}) reveals that if $s\in (0,s_\star)$ is such that
  $s\le s_1$, then 
  \bas
	\frac{w(s,t_0)}{\ow(s,t_0)} 
	\le \frac{c_1 s}{n} \cdot \frac{s^\frac{2}{n}+b_0}{2s}
	\le \frac{c_1}{2n} s_1^\frac{2}{n} + \frac{1}{2}
	= 1,
  \eas
  whereas if $s_1<s<s_\star$, then due to (\ref{22.4}) and (\ref{22.5}),
  \bas
	\ow(s,t_0)-w(s,t_0)
	&\ge& \frac{2s}{s^\frac{2}{n}+b_0} - \big(2s^{1-\frac{2}{n}} - c_2 \big) 
	\ge -2b_0 \max\{s_\star^{1-\frac{4}{n}},s_1^{1-\frac{4}{n}}\} + c_2 
	\ge 0.
  \eas
  Consequently,
  \bas
	\ow(s,t_0) \ge w(s,t_0)
	\qquad \mbox{for all } s\in (0,s_\star),
  \eas
  which together with (\ref{22.11}), (\ref{22.12}) and (\ref{22.10}) enables us to employ a comparison principle for
  (\ref{0w}) (\cite{win_IJM}) to infer that
  \bas
	\ow(s,t) \ge w(s,t)
	\qquad \mbox{for all $s\in (0,s_\star)$ and } t>t_0.
  \eas
  In light of (\ref{22.7}), this establishes (\ref{22.1}).
\qed
Through a Bernstein-type reasoning, in its basic strategy inspired by a precedent in \cite{win_MATANN}
operating on bounded domains, the pointwise bounds for the cumulated densities $w$ implied by Lemma \ref{lem22} and (\ref{pb})
can be seen to entail $L^\infty$ estimates for $u=nw_s$. This allows to prove 
%
our main result on the existence of a bounded set in $L^\infty(\R^n)$ that 
absorbs all trajectories emanating from initial data as in (\ref{init}).\abs
\proofc of Theorem \ref{theo25}. \quad
  Since $u(r,t)=nw_s(r^n,t)$ for $r>0$ and $t>0$, the statement immediately results if we prove the following\\
{\cred
\underline{Claim:} \
{\it  If $3\le n\le 9$, then there exists $C>0$ with the property that whenever $u_0$ satisfies (\ref{init}) with $u_0\not\equiv 0$,
  it is possible to find
  $t_2=t_2(u_0)>0$ such that
  \be{24.1}
	w_s(s,t) \le C
	\qquad \mbox{for all $s\ge 0$ and } t\ge t_2.
  \ee
} 
\noindent
  To achieve this, we
}
  take $B=B(1)>0$ be as provided by Lemma \ref{lem22}, and we claim that the intended conclusion holds if we let
  \be{24.2}
	C:=\sqrt{\frac{8}{n^2} \cdot \max \{c_1 \, , \, c_2\} \cdot \max \Big\{ \frac{2}{B} \, , \, 2\Big\}}
  \ee
  with
  \be{24.22}
	c_1:=\frac{14 n^2 (1+4\K^2)}{B} + \frac{4n^2(2-\frac{2}{n})^2}{B} + \frac{64}{B^3} + \frac{4n^2\K (3+\K)}{B} 
	+ \frac{4n(1+2\K)}{B^2}
	+ \frac{2\K}{B}
  \ee
  and
  \be{24.23}
	c_2:=14 n^2\big(1+2^{4-\frac{2}{n}} \K^2\big) + 4n^2\Big(2-\frac{2}{n}\Big)^2 + 64 
	+ 2n^2\Big( 2-\frac{2}{n}\Big) (1+4\K)
	+ 2^{4-\frac{2}{n}} n^2 \K (2+\K) +2\K,
  \ee
  where, again, $\K$ is as in (\ref{K}).\abs
  To confirm this, given $u_0$ with the properties in question we can employ Lemma \ref{lem22} to find $t_1=t_1(u_0)>0$ such that
  \bas
	w(s,t) \le \frac{2s}{s^\frac{2}{n}+B}
	\qquad \mbox{for all $s\in [0,1]$ and } t>t_1,
  \eas
  whence in particular,
  \be{24.3}
	w(s,t) \le \frac{2}{B} \cdot s
	\qquad \mbox{for all $s\in [0,1]$ and } t>t_1;
  \ee
  from (\ref{pb}) we moreover know that
  \be{24.4}
	w(s,t) \le 2s^{1-\frac{2}{n}}
	\qquad \mbox{for all $s>1$ and } t>t_1.
  \ee
  With $\chi$ and $(\ze)_{\eps\in (0,1)}$ as in (\ref{chi}) and (\ref{ze}), we now let
  \be{24.44}
	\xe(s,t):=\chi(t-t_1) \cdot s \ze^2(s),
	\qquad s\ge 0, \ \eps\in (0,1),
  \ee
  and
	\[
	\zeps(s,t):=\xe(s,t)\cdot \frac{w_s^2(s,t)}{w(s,t)},
	\qquad s>0, \ t\ge t_1, \ \eps\in (0,1),
	\]
  noting that then $\zeps$ is well-defined for any such $\eps$ due to the positivity of $w$ in $(0,\infty)^2$ resulting
  from the strong maximum principle, and that moreover
  \bas
	\supp \zeps \subset \Big[\frac{\eps}{2},\frac{2}{\eps}\Big] \times [t_1+1,\infty)
	\qquad \mbox{for all } \eps\in (0,1)
  \eas
  according to (\ref{chi}) and (\ref{z2}).
  Consequently, for each $T>t_1+1$ we can find $\wh{s}\in [\frac{\eps}{2},\frac{2}{\eps}]$ and $\wh{t}\in [t_1+1,T]$ such that the 
  positive number $S:=\|\zeps\|_{L^\infty((0,\infty)\times (t_1,T))}$ satisfies
	\[
	S=\zeps(\wh{s},\wh{t}),
	\]
  and that $z_{\eps t}(\wh{s},\wh{t})\ge 0, z_{\eps s}(\wh{s},\wh{t})=0$ and $z_{\eps ss}(\wh{s},\wh{t})\le 0$.
  Since
  \be{24.7}
	z_{\eps s} = 2\xe w^{-1} w_s w_{ss} - \xe w^{-2} w_s^3 + \xes w^{-1} w_s^2
  \ee
  and thus
	\bas
	z_{\eps ss}
	&=& 2\xe w^{-1} w_s w_{sss}
	+ 2\xe w^{-1} w_{ss}^2
	- 5\xe w^{-2} w_s^2 w_{ss}
	+ 4\xes w^{-1} w_s w_{ss} \nn\\
	& & + 2\xe w^{-3} w_s^4 
	- 2\xes w^{-2} w_s^3
	+ \xess w^{-1} w_s^2
	\eas
  for $s>0$ and $t>t_1$, in view of (\ref{0w}) this implies that at $(s,t)=(\wh{s},\wh{t})$,
  \bea{24.9}
	0\le z_{\eps t}
	&=& 2\xe w^{-1} w_s \cdot \Big\{ n^2 s^{2-\frac{2}{n}} w_{sss} + n^2\Big(2-\frac{2}{n}\Big) s^{1-\frac{2}{n}} w_{ss}
		+ nww_{ss} + nw_s^2 \Big\} \nn\\
	& & - \xe w^{-2} w_s^2 \cdot \Big\{ n^2 s^{2-\frac{2}{n}} w_{ss} + nww_s \Big\} \nn\\
	& & + \xet w^{-1} w_s^2 \nn\\
	&\le& - 2n^2 \xe s^{2-\frac{2}{n}} w^{-1} w_{ss}^2 + 4n^2 \xe s^{2-\frac{2}{n}} w^{-2} w_s^2 w_{ss} \nn\\
	& & - 4n^2 \xes s^{2-\frac{2}{n}} w^{-1} w_s w_{ss}
	+ 2n^2\Big(2-\frac{2}{n}\Big) \xe s^{1-\frac{2}{n}} w^{-1} w_s w_{ss} 
	+ 2n\xe w_s w_{ss} 
	\nn\\
	& & -2n^2 \xe s^{2-\frac{2}{n}} w^{-3} w_s^4 \nn\\
	& & + 2n^2 \xes s^{2-\frac{2}{n}} w^{-2} w_s^3
	+ n\xe w^{-1} w_s^3 \nn\\
	& & - n^2 \xess s^{2-\frac{2}{n}} w^{-1} w_s^2 + \xet w^{-1} w_s^2.
  \eea
  But in view of (\ref{24.7}), at each point where $\xe\ne 0$ the identity $z_{\eps s}=0$ is equivalent to the relation
  \bas
	w_{ss}=\frac{1}{2} w^{-1} w_s^2 - \frac{\xes}{2\xe} w_s,
  \eas
  so that as $(s,t)=(\wh{s},\wh{t})$, (\ref{24.9}) entails that
  \bas
	0
	&\le& -2n^2 \xe s^{2-\frac{2}{n}} w^{-1} \cdot \Big\{ \frac{1}{2} w^{-1} w_s^2 - \frac{\xes}{2\xe} w_s \Big\}^2 \nn\\
	& & + 4n^2 \xe s^{2-\frac{2}{n}} w^{-2} w_s^2 \cdot \Big\{ \frac{1}{2} w^{-1} w_s^2 - \frac{\xes}{2\xe} w_s\Big\} \nn\\
	& & - 4n^2 \xes s^{2-\frac{2}{n}} w^{-1} w_s \cdot \Big\{ \frac{1}{2} w^{-1} w_s^2 - \frac{\xes}{2\xe} w_s\Big\} \nn\\
	& & + 2n^2 \Big(2-\frac{2}{n}\Big) \xe s^{1-\frac{2}{n}} w^{-1} w_s \cdot 
		\Big\{ \frac{1}{2} w^{-1} w_s^2 - \frac{\xes}{2\xe} w_s\Big\} \nn\\
	& & 
	+ 2n \xe w_s \cdot \Big\{ \frac{1}{2} w^{-1} w_s^2 - \frac{\xes}{2\xe} w_s\Big\} 
	-2n^2 \xe s^{2-\frac{2}{n}} w^{-3} w_s^4 
	\nn\\
	& & 
	+ 2n^2 \xes s^{2-\frac{2}{n}} w^{-2} w_s^3
	+ n\xe w^{-1} w_s^3 
	- n^2 \xess s^{2-\frac{2}{n}} w^{-1} w_s^2 + \xet w^{-1} w_s^2 \nn\\
	&=& - \frac{n^2}{2} \xe s^{2-\frac{2}{n}} w^{-3} w_s^4 \nn\\
	& & - n^2 \xes s^{2-\frac{2}{n}} w^{-2} w_s^3
	+ n^2\Big(2-\frac{2}{n}\Big) \xe s^{1-\frac{2}{n}} w^{-2} w_s^3
	+ 2 n\xe w^{-1} w_s^3 
	\nn\\
	& & + \frac{3}{2} n^2 \frac{\xes^2}{\xe} s^{2-\frac{2}{n}} w^{-1} w_s^2
	- n^2\Big(2-\frac{2}{n}\Big) \xes s^{1-\frac{2}{n}} w^{-1} w_s^2
	- n^2 \xess s^{2-\frac{2}{n}} w^{-1} w_s^2 \nn\\
	& & 
	- n\xes w_s^2
	+ \xet w^{-1} w_s^2.
  \eas
  Upon multiplication by $\frac{s^{-2+\frac{2}{n}} w^2}{w_s^2}$, this shows that at this point, due to Young's inequality 
  we have
  \bas
	\frac{n^2}{2} \xe w^{-1} w_s^2
	&\le& - n^2 \xes w_s
	+ n^2\Big(2-\frac{2}{n}\Big) \xe s^{-1} w_s 
	+ 2n\xe s^{-2+\frac{2}{n}} ww_s 
	\nn\\
	& & + \frac{3}{2}n^2 \frac{\xes^2}{\xe} w
	- n^2\Big(2-\frac{2}{n}\Big) \xes s^{-1} w
	- n^2 \xess w  
	- n \xes s^{-2+\frac{2}{n}} w^2 
	+ \xet s^{-2+\frac{2}{n}} w  \nn\\
	&\le& \frac{n^2}{8} \xe w^{-1} w_s^2
	+ 2n^2 \frac{\xes^2}{\xe} w 
	+ \frac{n^2}{8} \xe w^{-1} w_s^2
	+ 2n^2\Big(2-\frac{2}{n}\Big)^2 \xe s^{-2} w \nn\\
	& & +\frac{n^2}{8} \xe w^{-1} w_s^2
	+ 8\xe s^{-4+\frac{4}{n}} w^3 
	+ \frac{3}{2} n^2 \frac{\xes^2}{\xe} w 
	- n^2\Big(2-\frac{2}{n}\Big) \xes s^{-1} w 
	- n^2 \xess w \nn\\
	& & 
	- n \xes s^{-2+\frac{2}{n}} w^2
	+  \xet s^{-2+\frac{2}{n}} w 
  \eas
  and thus
  \bea{24.10}
	\frac{n^2}{8} \xe w^{-1} w_s^2
	&\le& \frac{7}{2} n^2 \frac{\xes^2}{\xe} w
	+ 2n^2\Big(2-\frac{2}{n}\Big)^2 \xe s^{-2} w
	+ 8\xe s^{-4+\frac{4}{n}} w^3 
	\nn\\
	& & - n^2\Big(2-\frac{2}{n}\Big) \xes s^{-1} w
	- n^2 \xess w \nn\\
	& & 
	- n \xes s^{-2+\frac{2}{n}} w^2
	+ \xet s^{-2+\frac{2}{n}} w.
  \eea
  We now make use of (\ref{24.44}) and (\ref{chi1}) to see that
  \bas
	0\le \xe(s,t) \le s
  \eas
  and
  \bas
	\xes(s,t) =\chi(t-t_1) \ze^2(s) + 2\chi(t-t_1) s\ze(s)\zes(s)
  \eas
  as well as
  \bas
	\xess(s,t)=4\chi(t-t_1) \ze(s) \zes(s) + 2\chi(t-t_1) s \ze(s) \zess(s) + 2\chi(t-t_1) s\zes^2(s)
  \eas
  and
  \bas
	|\xet(s,t)|
	= |\chi'(t-t_1) s\ze^2(s)| \le \K s
  \eas
  for $s>0, t>t_1$ and $\eps\in (0,1)$, whence in particular, by Young's inequality an the fact that $0\le \chi\le 1$,
	\bas
	\frac{\xes^2(s,t)}{\xe(s,t)}
	&\le& \frac{2\cdot \Big\{ \chi^2(t-t_1) \ze^4(s) + 4\chi^2(t-t_1) s^2 \ze^2(s) \zes^2(s)\Big\}}{\chi(t-t_1) s\ze^2(s)} \nn\\
	&\le& 2\frac{\ze^2(s)}{s} + 8 s\zes^2(s) \qquad  \mbox{whenever $\eps\in (0,1)$, $s>0$ and $t>t_1$ are such that $\xe(s,t)>0$,}
	\eas
  and 
  \be{24.111}
	|\xes(s,t)| \le \ze^2(s) + 2s\ze(s) |\zes(s)|
	\qquad \mbox{for all $s>0, t>t_1$ and } \eps\in (0,1)
  \ee
  as well as
  \[
	|\xess(s,t)| \le 4|\zes(s)| + 2s|\zess(s)| + 2s\zes^2(s) 
	\qquad \mbox{for all $s>0, t>t_1$ and } \eps\in (0,1).
  \]
  In the case when $\wh{s}\le 1$, on the right-hand side of (\ref{24.10}) we can therefore conclude, using (\ref{z3}) along with  
  (\ref{24.3}), the inequality $2-\frac{2}{n}\ge 0$ and the nonnegativity of $\zes$ on $(0,1)$, that
  \bas
	\frac{7}{2}n^2 \frac{\xes^2(\whs,\wht)}{\xe(\whs,\wht)} w(\whs,\wht)
	&\le& \frac{7}{2} n^2 \cdot \Big\{ 2\frac{\ze^2(\whs)}{\whs} + 8\whs \zes^2(\whs)\Big\} \cdot \frac{2}{B} \whs \nn\\
	&\le& \frac{7 n^2}{B} \cdot \Big\{ 2 + 8 \eps^2 \cdot \Big(\frac{\K}{\eps} \Big)^2\Big\} \\
	&=& \frac{14 n^2 (1+4 \K^2)}{B}
  \eas
  and
  \bas
	2n^2\Big(2-\frac{2}{n}\Big)^2 \xe(\whs,\wht) \whs^{-2} w(\whs,\wht)
	\le 2n^2 \Big(2-\frac{2}{n}\Big)^2 \cdot \whs \cdot \whs^{-2} \cdot \frac{2}{B} \whs
	= \frac{4n^2 (2-\frac{2}{n})^2}{B}
  \eas
  and
  \bas
	8\xe(\whs,\wht) \whs^{-4+\frac{4}{n}} w^3( \whs,\wht)
	\le 8\whs\cdot\whs^{-4+\frac{4}{n}} \cdot \Big(\frac{2}{B} \whs\Big)^3
	= \frac{64}{B^3} \whs^\frac{4}{n}
	\le \frac{64}{B^3}
  \eas
  and
  \bas
	- n^2\Big(2-\frac{2}{n}\Big) \xes(\whs,\wht) \whs^{-1} w(\whs,\wht) 
	\le - 2n^2\Big(2-\frac{2}{n}\Big)\chi(t-t_1)\ze(\whs)\zes(\whs)\whs^{-1} w(\whs,\wht)
	\le 0
  \eas
  as well as
  \bas
	- n^2 \xess(\whs,\wht) w(\whs,\wht)
	&\le& n^2 \cdot \Big\{ 4|\zes(\whs)| + 2\whs |\zess(\whs)| + 2\whs \zes^2(\whs)\Big\} \cdot \frac{2}{B} \whs \nn\\
	&\le& \frac{2n^2}{B} \cdot \Big\{ 4\eps \cdot \frac{\K}{\eps} + 2\eps^2 \cdot \frac{\K}{\eps^2}
		+ 2\eps^2\cdot\Big(\frac{\K}{\eps}\Big)^2 \Big\} \\
	&=& \frac{4n^2 \K (3+\K)}{B}
  \eas
  and
  \bas
	- n\xes \whs^{-2+\frac{2}{n}} w^2(\whs,\wht)
	&\le& n \cdot \Big\{ \ze^2(\whs) + 2\whs \ze(\whs) |\zes(\whs)| \Big\} \cdot \whs^{-2+\frac{2}{n}} w^2(\whs,\wht) \\
	&\le& n\cdot \Big\{1 + 2\eps \cdot \frac{\K}{\eps}\Big\} \cdot \whs^{-2+\frac{2}{n}} \cdot \Big(\frac{2}{B} \whs\Big)^2 \\
	&\le& \frac{4n(1+2\K)}{B^2}
  \eas
  and
  \bas
	\xet(\whs,\wht) \whs^{-2+\frac{2}{n}} w(\whs,\wht)
	\le \K \whs \cdot \whs^{-2+\frac{2}{n}} \cdot \frac{2}{B} \whs
	= \frac{2\K}{B} \whs^\frac{2}{n} 
	\le \frac{2\K}{B},
  \eas
  so that according to (\ref{24.22}) and \eqref{24.10}
  \be{24.13}
	\frac{n^2}{8} \cdot \xe(\whs,\wht) w^{-1}(\whs,\wht) w_s^2(\whs,\wht) \le c_1
	\qquad \mbox{if } \whs \le 1.
  \ee
  In the case when, alternatively, $\whs> 1$, we instead rely on (\ref{z4}) and (\ref{24.4}), and additionally
  on (\ref{24.111}), to similarly see that
  \bas
	\frac{7}{2} n^2 \frac{\xes^2(\whs,\wht)}{\xe(\whs,\wht)} w(\whs,\wht)
	&\le& \frac{7}{2} n^2 \cdot \Big\{ 2\frac{\ze^2(\whs)}{\whs} + 8 \whs \zes^2(\whs)\Big\} \cdot 2\whs^{1-\frac{2}{n}} \\
	&\le& 14 n^2 \cdot \Big\{ 1+4\cdot\Big(\frac{2}{\eps}\Big)^{2-\frac{2}{n}} \cdot (\K\eps)^2 \Big\} \\
	&\le& 14 n^2 (1+2^{4-\frac{2}{n}} \K^2)
  \eas
  and
  \bas
	2n^2\Big(2-\frac{2}{n}\Big)^2 \xe(\whs,\wht) \whs^{-2} w(\whs,\wht)
	\le 2n^2 \Big(2-\frac{2}{n}\Big)^2 \cdot \whs \cdot \whs^{-2} \cdot 2\whs^{1-\frac{2}{n}} 
	\le 4n^2\Big(2-\frac{2}{n}\Big)^2,
  \eas
  that
  \bas
	8\xe(\whs,\wht) \whs^{-4+\frac{4}{n}} w^3( \whs,\wht)
	\le 8\whs\cdot\whs^{-4+\frac{4}{n}} \cdot (2\whs^{1-\frac{2}{n}})^3
	= 64\whs^{-\frac{2}{n}} \le 64
  \eas
  and
  \bas
	- n^2\Big(2-\frac{2}{n}\Big) \xes(\whs,\wht) \whs^{-1} w(\whs,\wht) 
	&\le& n^2\Big(2-\frac{2}{n}\Big) \cdot \Big\{ \ze^2(\whs) + 2\whs \ze(\whs) |\zes(\whs)| \Big\} \cdot \whs^{-1} w(\whs,\wht)
		\\
	&\le& n^2\Big(2-\frac{2}{n}\Big) \cdot \Big\{ 1 + 2\cdot\frac{2}{\eps} \cdot \K \eps \Big\} \cdot \whs^{-1} \cdot 
		2 \whs^{1-\frac{2}{n}} \\
	&\le& 2n^2\Big(2-\frac{2}{n}\Big) \cdot (1+4\K)
  \eas
  and
  \bas
	- n^2 \xess(\whs,\wht) w(\whs,\wht)
	&\le& n^2 \cdot \Big\{ 4|\zes(\whs)| + 2\whs |\zess(\whs)| + 2\whs \zes^2(\whs)\Big\} \cdot 2\whs^{1-\frac{2}{n}} \\
	&\le& 2n^2 \cdot \Big\{ 4\cdot \Big(\frac{2}{\eps}\Big)^{1-\frac{2}{n}} \cdot \K \eps
	+ 2\cdot\Big( \frac{2}{\eps}\Big)^{2-\frac{2}{n}} \cdot \K \eps^2
	+ 2\cdot\Big(\frac{2}{\eps}\Big)^{2-\frac{2}{n}} \cdot (\K\eps)^2 \Big\} \\
	&\le& 2^{4-\frac{2}{n}} n^2 \K (2+\K),
  \eas
  and that
  \bas
	\xet(\whs,\wht) \whs^{-2+\frac{2}{n}} w(\whs,\wht)
	\le \K \whs \cdot \whs^{-2+\frac{2}{n}} \cdot 2\whs^{1-\frac{2}{n}}
	= 2\K.
  \eas
  In line with (\ref{24.23}), from (\ref{24.10}) we thus infer that
  \bas
	\frac{n^2}{8} \cdot \xe(\whs,\wht) w^{-1}(\whs,\wht) w_s^2(\whs,\wht) \le c_2
	\qquad \mbox{if } \whs > 1, 
  \eas
  which combined with (\ref{24.13}) shows that
  \bas
	\chi(t-t_1) s\ze^2(s) \cdot \frac{w_s^2(s,t)}{w(s,t)}
	\le c_3:=\frac{8}{n^2} \cdot \max \{ c_1 \, , \, c_2\}
	\qquad \mbox{whenever $\eps\in (0,1), s>0$ and $t>t_1$},
  \eas
  because $T>t_1$ was arbitrary. 
  Since $c_3$ does not depend on $\eps\in (0,1)$, and since $\chi\equiv 1$ on $[2,\infty)$ by (\ref{chi}), this entails that
  \bas
	w_s^2(s,t) \le c_3 \cdot \frac{w(s,t)}{s}
	\qquad \mbox{for all $s>0$ and } t>t_1+2,
  \eas
  so that a final application of (\ref{24.3}) and (\ref{24.4}) reveals that
  \bas
	w_s^2(s,t) \le \frac{2c_3}{B}
	\qquad \mbox{for all $s\in (0,1]$ and } t>t_1+2,
  \eas
  and that
  \bas
	w_s^2(s,t) \le 2c_3 s^{-\frac{2}{n}} \le 2c_3
	\qquad \mbox{for all $s> 1$ and } t>t_1+2.
  \eas
  By definition of $c_3$, these inequalities precisely assert (\ref{24.1}) with $C$ as in (\ref{24.2}), and with 
  $t_2\equiv t_2(u_0):=t_1+2$.
\qed
%
%
%
%
%
%
%
%
The intended explicit instability property has actually been revealed by Lemma \ref{lem22}:\abs
\proofc of Proposition \ref{prop26}. \quad
  We let $s_0:=r_0^n$ and take $B=B(s_0)>0$ as accordingly provided by Lemma \ref{lem22}.
  Then given any $u_0$ fulfilling (\ref{init}), from Lemma \ref{lem22} we obtain $t_1=t_1(u_0)>0$ such that 
  $w(s,t) \le \frac{2s}{s^\frac{2}{n}+B}$ for all $s\in (0,s_0)$ and $t\ge t_1$, which in view of (\ref{w}) means that
  \bas
	\Mint_{B_r(0)} u(x,t) dx 
	= \frac{n}{r^n} w(r^n,t) 
	\le \frac{n}{r^n} \cdot \frac{2r^n}{r^2+B}
	= \frac{2n}{r^2+B}
	\qquad \mbox{for all $r\in (0,r_0)$ and } t\ge t_1.
  \eas
  Since, as can readily be verified,
  \bas
	\Mint_{B_r(0)} u_\star(x) dx = \frac{2n}{r^2}
	\qquad \mbox{for all } r>0,
  \eas
  this implies that
  \bas
	\Mint_{B_r(0)} u_\star(x) dx - \Mint_{B_r(0)} u(x,t) dx
	&\ge& 2n\cdot \Big(\frac{1}{r^2}-\frac{1}{r^2+B}\Big) \\
	&=& \frac{2nB}{r^2(r^2+B)}
	\ge \frac{2nB}{r_0^2 (r_0^2+B)}
	\qquad \mbox{for all $r\in (0,r_0)$ and } t\ge t_1,
  \eas
  and thereby yields the claim.
\qed

\bigskip

{\bf Acknowledgement.} \quad
The second author 
acknowledges support of the {\em Deutsche Forschungsgemeinschaft} (Project No.~462888149) and of the Department of Mathematics of the University of Bologna for his visit to Bologna where parts of this work have been achieved.
The first named author was partially supported by the INdAM-GNAMPA Project 2022 {\em Studi asintotici in problemi parabolici ed ellittici}.

\end{document}